\documentclass[reqno]{amsart}
\usepackage{calrsfs} 
\usepackage{cite} 
\usepackage{graphicx} 
\usepackage{indentfirst} 
\usepackage{xcolor}

\newtheorem{ntheorem}{Theorem}
\theoremstyle{definition}
\newtheorem{theorem}{Theorem}
\renewcommand*{\thetheorem}{\Alph{theorem}}
\newtheorem{lemma}[theorem]{Lemma}
\newtheorem{corollary}{Corollary}

\newtheorem*{remark}{Remark}

\usepackage{caption}
\theoremstyle{definition}

\makeatletter 
\def\thmhead@plain#1#2#3{%
  \thmname{#1}\thmnumber{\@ifnotempty{#1}{ }\@upn{#2}}%
  \thmnote{{\the\thm@notefont#3}}}
\let\thmhead\thmhead@plain
\makeatother

\begin{document}
\title[Extremal polynomials for the Rogosinski--Szeg\H{o} problem]
{Extremal polynomials for the Rogosinski--Szeg\H{o} estimates of the 
third coefficient of nonnegative sine polynomials}

\author[D.~Dmitrishin]{Dmitriy Dmitrishin}
\address{Dmitriy Dmitrishin: Odessa National Polytechnic University, 1 Shevchenko Ave., Odesa 65044, Ukraine}
\email{dmitrishin@op.edu.ua}

\author[A.~Stokolos]{Alexander Stokolos}
\address{Alexander Stokolos: Georgia Southern University, Statesboro GA, 30460, USA }
\email{astokolos@georgiasouthern.edu}

\author[W.~Trebels]{Walter Trebels}
\address{Walter Trebels: Department of Mathematics, AG Algebra, Technical University of Darmstadt, 64289 Darmstadt, Germany}
\email{trebels@mathematik.tu-darmstadt.de}

\date{\today}

\begin{abstract}
In the class of  normalized sine-polynomials $S(t),$ non-negative on $[0,\pi],$ 
W.~Rogosinski and G.~Szeg\H{o} 1950 considered a number of extremal problems and proved, among other things, sharp upper and lower estimates for the coefficient $a_3.$ Their proof is based on the Luk\'acs representation of non-negative algebraic polynomials. This method does not lead to the construction of polynomials attaining the extreme values.

We consider the corresponding problem in the framework of normalized typically real polynomials $P(z)$ on the unit disc in $\mathbb C.$
By L.~Fej\'er's method with the additional use of
the Chebyshev polynomials of the second kind and their derivatives we regain the sharp upper and lower estimates for $a_3$ and identify the extremal polynomials. The corresponding statements for sine polynomials follow by the observation $S(t)=\text{Im}\{P(e^{it})\}$. For odd $N$ the extremizers are unique, for even $N$ there is a one-parameter family of extremizers. 
\end{abstract}
\maketitle


\section{Introduction}
In the class of normalized sine-polynomials, non-negative on $[0,\pi],$ 
\begin{equation}\label{eq:sine}
S(t)=\sin{t}+\sum_{j=2}^N a_j \sin{jt},\quad a_N\ne0,
\end{equation}
W.~Rogosinski and G.~Szeg\H{o} \cite{RS50} considered a number of extremal problems and 
proved, among other things, that 
\begin{align}\label{eq:bounds}
\begin{split}
1+2\cos\frac{\pi}{n'+3}&\le a_3\le1+2\cos\frac{2\pi}{n'+3}
\quad\text{if $n'$ is even},\\
1-2\cos\vartheta&\le a_3\le1+2\cos\frac{2\pi}{n'+3}
\quad\text{if $n'$ is odd},  
\end{split}
\end{align}
where $n'=[(N-1)/2]$, the function $[x]$ defines the integer part of a number $x$,
and $\vartheta$ is the minimum positive root of the equation 
\[
(n'+4)\cos\frac{(n'+2)\vartheta}{2}+(n'+2)\cos\frac{(n'+4)\vartheta}{2}=0.
\]
The proof of these estimates is based on the Luk\'acs representation of 
non-negative algebraic polynomials and a brilliant application of fine
variational methods of mechanical quadratures. However, the problem of constructing extremizers (sine polynomials which attain the bounds in \eqref{eq:bounds}) was not solved. But Rogosinski and Szeg\H{o} \cite{RS50} discussed possible ways to obtain these. They noted that the approach via the Fej\'er--Riesz representation of a non-negative
trigonometric polynomial~\cite[6.5, Problem~41]{PS98} and a subsequent application 
of the Rayleigh method "is not easily adaptable for
obtaining explicit results, in particular when degree $N$ is large.'' 

In \cite{DST24} we consider the problem in the framework of geometric complex analysis. The set of univalent functions in the unit disc $\mathbb{D}=\{z\in\mathbb{C}:|z|<1\}$ is not a linear space. To bypass this inconvenience, W.~Rogosinski in~\cite{WR32} introduced a class $\mathcal T$ of typically real functions.  $\mathcal T$ consists of all holomorphic functions  $F$ on $\mathbb D$ which satisfy

i) $\forall z\in \mathbb R\cap\mathbb D,\; F(z)\in \mathbb R,$ 

ii) $\forall z\in \mathbb{D} \backslash\mathbb R,\; \text{Im}\{F(z)\}\cdot\text{Im}\{z\}>0.$ 

The class of typically real functions is convex, which appears to be useful when
solving various extremal problems. Moreover, many extremal estimates in the class
of typically real functions remain the same or are similar to the corresponding
estimates in the class of univalent functions.

Denote by $\mathcal{T}_N$ the set of normalized typically real polynomials of degree $N$ 
$$P(z)=z+\sum_{j=2}^N a_j z^j,\quad a_N\not=0, \quad \mbox{on $\mathbb{D}$}.$$
Obviously, all the coefficients of a typically real polynomial are real and $\text{Im}\{P(e^{it})\}$ is just the sine polynomial \eqref{eq:sine}.

In \cite{DST24} it is shown by a bigger technical effort that this method works for the second coefficient $a_2$ to deduce the (known) sharp bounds, to obtain a construction of the coefficients of the extremizers as well as to represent the extremizers in the form of rational functions without poles. Here it turns out that the method we used for $a_2$ also works for $a_3.$ 
We have to pay for this with even higher technical efforts. Some proofs contain cumbersome trigonometric transformations, which are omitted in the paper. Interested reader can validate formulas, e.g. the core formula \eqref{eq:det}, using MAPLE of MATHEMATICA software.

The interplay between the typically real extremizers and the sine extremizers is discussed in Sec. \ref{sec:5}.

The asymmetry of the lower and upper bounds in \eqref{eq:bounds} let us expect that we have to deal with more than the two cases occurring in the discussion of the second coefficient $a_2$ (see \cite{DST24}). For clarification, we separate the hidden cases in the parameter $n'$ of \eqref{eq:bounds}.

$$
n'=\left[\frac{N-1}2\right] \textit{is odd}\qquad \Longleftrightarrow\qquad 
\left\{
\begin{array}{l}
N\textit{ is odd},\ (N+3)/2\textit{ is odd}\\
N\textit{ is even},\ (N+2)/2\textit{ is odd},
\end{array}
\right.
$$

$$
n'=\left[\frac{N-1}2\right] \textit{is even}\qquad \Longleftrightarrow\qquad 
\left\{
\begin{array}{l}
N\textit{ is odd},\ (N+3)/2\textit{ is even}\\
N\textit{ is even},\ (N+2)/2\textit{ is even}.
\end{array}
\right.
$$
Hence, on account of our approach, we need to consider the four cases listed on the right-hand side of the two preceding displays. This implies that we have to construct four extremizers for the upper bounds of $a_3;$ we mark the corresponding problems by  $A_\text{max}$, $B_\text{max}$,
$C_\text{max}$, $D_\text{max}.$ Analogously, we denote the corresponding problems for the lower bound by  $A_\text{min}$, $B_\text{min}$, $C_\text{min}$, $D_\text{min}.$ 

It turns out that the cases $A_\text{max}$, $B_\text{max}$ as well as $C_\text{max}$, $D_\text{max}$ coincide. 
{ Thus, we have to consider the following problem:\\

\begin{minipage}{5cm}
\vspace{-.5cm}
\hspace{0.6cm} Find the \textit{maxima} for
\begin{enumerate}
\item[a)] $A_\text{max}$, $B_\text{max}$, $N$ is odd;
\item[b)] $C_\text{max}$, $D_\text{max}$, $N$ is even;\\
\end{enumerate}
\end{minipage}
\begin{minipage}{8cm}
\hspace{0.6cm} Find the \textit{minima} for
\begin{enumerate}
\item[c)] $A_\text{min}$, $N$ is odd, $(N+3)/2$ is odd;
\item[d)] $B_\text{min}$, $N$ is odd, $(N+3)/2$ is even; 
\item[e)] $C_\text{min}$, $N$ is even, $(N+2)/2$ is odd;
\item[f)] $D_\text{min}$, $N$ is even, $(N+2)/2$ is even.
\end{enumerate}
\end{minipage}
}
The results of this paper are to be seen in the context of previous contributions by several authors to this area of research - in this regard, see the Introduction of \cite{DST24}.

\section{Main results}
\subsection{Coefficients of the extremal typically real polynomials}
For their construction the Chebyshev polynomials of the second kind 
\[
U_j(x)=U_j(\cos{t})=\frac{\sin(j+1)t}{\sin{t}}=2^j x^j+\ldots
\]
and their derivatives $U'_j(x)$ play an essential role. Define for odd indices $(k\in\mathbb N_0)$
\begin{align*}
A_{2k+1}&(y)=U_{2k}(y)+\alpha_N\cdot\frac{y}{1-y^2}
\left(
y+U_{2k-1}(y)\frac{3y^2-1}{2y^2}-U_{2k}(y)\frac{k(1-y^2)+y^2}{y}
\right),
\end{align*}
where { $U_{-1}(x)\equiv 0$} and
\[
\alpha_N=\left\{
\begin{array}{l}
2/(N+5)\textit{ if $N$ is odd},\\
2/(N+4)\textit{ if $N$ is even};
\end{array}
\right.
\]
define for even $N$ and even indices
\[
A_{2k+2}(y)=\frac4{N+4}\cdot y(1-y^2)U_k(y)U'_{N/2-k}(y),
\quad k=0,\ldots,\frac{N-2}2.
\]
\begin{theorem}\label{th:A}
\hspace{1cm}
{}
\\
\hspace*{.4cm}\textit{a)} Cases $A_\text{max}$, $B_\text{max}$: Let \textit{$N$ be odd} and
$y_N^{(1)}$ be the maximum root of the equation $U_{(N+3)/2}(x)=0$. Then 
\[
\max_{P\in\mathcal{T}_N}\{a_3\}=4(y_N^{(1)})^2-1,\qquad y_N^{(1)}=\cos\frac{2\pi}{N+5}.
\] 
The extremal polynomial $\displaystyle P(z)=z+\sum_{j=2}^N a_j z^j$ is unique, its coefficients 
are defined by the formulas 
$$a_{2k+1}=A_{2k+1}(y_N^{(1)}),\qquad a_{2k+2}=0,\quad k=0,\ldots,(N-1)/2.$$

\textit{b)} Cases $C_\text{max}$, $D_\text{max}$: Let \textit{$N$ be even} and
$y_N^{(3)}$ be the maximum root of the equation $U_{(N+2)/2}(x)=0.$ Then 
\[
\max_{P\in\mathcal{T}_N}\{a_3\}=4(y_N^{(3)})^2-1,\qquad y_N^{(3)}=\cos\frac{2\pi}{N+4}.
\]
The collection of the extremal polynomials $\displaystyle P(z)=z+\sum_{j=2}^N a_j z^j$ forms a
one-parameter family, the coefficients of these polynomials are defined by the
formulas 
$$a_{2k+1}=A_{2k+1}(y_N^{(3)}),\; a_{2k+2}=\tau A_{2k+2}(y_N^{(3)}),\;
\tau\in[-1,1], \quad k=0,\ldots,(N-2)/2.$$

\textit{c)} Case $A_\text{min}$: Let \textit{$N$ be odd, $(N+3)/2$ be odd} and
$y_N^{(5)}$ be the minimum positive root of the equation $U'_{(N+3)/2}(x)=0.$
 Then
\[
\min_{P\in\mathcal{T}_N}\{a_3\}=4(y_N^{(5)})^2-1.
\]
The extremal polynomial $\displaystyle  P(z)=z+\sum_{j=2}^Na_jz^j$ is unique. 

For the construction of the coefficients introduce the auxiliary function 
\[
\zeta_k(x)=U_{k-1}(x)U_k(x)-U_{\frac{N+1}{2}-k}(x)U_{\frac{N+3}{2}-k}(x)
+\left(1-\frac{4k}{N+3}\right)U_\frac{N+1}{2}(x)U_\frac{N+3}{2}(x),
\]
$k=1,\ldots,(N+1)/2$ and define $z_k=\zeta_k(y_N^{(5)}).$ Then 
\[
a_{2k-1}=\frac
{\sum\limits_{j=1}^{\frac{N+3}{2}-k}z_j z_{j+k-1}
-\sum\limits_{j=1}^{\frac{N+1}{2}-k}z_j z_{j+k}}
{\sum\limits_{j=1}^{\frac{N+1}{2}}z_j z_j
-\sum\limits_{j=1}^{\frac{N-1}{2}}z_j z_{j+1}},\qquad
\ a_{2k}=0,
\ k=1,\ldots,\frac{N+1}{2}.
\]

\textit{d)} Case $B_\text{min}$: Let \textit{$N$ be odd, $(N+3)/2$ be even}  and
$y_N^{(2)}$ be the minimum positive root of the equation $U_{(N+3)/2}(x)=0.$  Then 
\[
\min_{P\in\mathcal{T}_N}\{a_3\}=4(y_N^{(2)})^2-1, \qquad y_N^{(2)}=\sin\frac{\pi}{N+5}.
\]
The extremal polynomial $\displaystyle P(z)=z+\sum_{j=2}^N a_j z^j$ is unique. Its coefficients 
are defined by the formulas 
$$a_{2k+1}=A_{2k+1}(y_N^{(2)}),\qquad a_{2k+2}=0,\quad k=0,\ldots,(N-1)/2.$$

\textit{e)} Case $C_\text{min}$: Let \textit{$N$ be even, $(N+2)/2$ be odd} and
$y_N^{(6)}$ be the minimum positive root of the equation $U'_{(N+2)/2}(x)=0.$ 
Then
\[
\min_{P\in\mathcal{T}_N}\{a_3\}=4(y_N^{(6)})^2-1.
\]
The set of the extremal polynomials $ P(z)=z+\sum_{j=2}^Na_jz^j$ forms a one-parameter family.
Analogously to \textit{c)} introduce for $k=1,2,...,N/2$
\[
\zeta^*_k(x)=U_{k-1}(x)U_k(x)-U_{\frac{N}{2}-k}(x)U_{\frac{N+2}{2}-k}(x)
+\left(1-\frac{4k}{N+2}\right)U_\frac{N}{2}(x)U_\frac{N+2}{2}(x),
\]
and define $z_k=\zeta^*_k(y_N^{(6)})$. Then,  for $\tau\in[-1,1]$ and $k=1,\ldots,N/2,$
\[
a_{2k-1}=\frac
{\sum\limits_{j=1}^{\frac{N+2}{2}-k}z_j z_{j+k-1}
-\sum\limits_{j=1}^{\frac{N}{2}-k}z_j z_{j+k}}
{\sum\limits_{j=1}^{\frac{N}{2}}z_j z_j
-\sum\limits_{j=1}^{\frac{N-2}{2}}z_j z_{j+1}},
\qquad 
a_{2k}=\tau\frac
{\sum\limits_{j=1}^{\frac{N+2}{2}-k}z_j z_{j+k-1}
-\sum\limits_{j=1}^{\frac{N-2}{2}-k}z_j z_{j+k+1}}
{\sum\limits_{j=1}^{\frac{N}{2}}z_j z_j
-\sum\limits_{j=1}^{\frac{N-2}{2}}z_j z_{j+1}}.
\]

\textit{f)} Case $D_\text{min}$: Let \textit{$N$ be even, $(N+2)/2$ be even} and
$y_N^{(4)}$ be the minimum positive root of the equation $U_{(N+2)/2}(x)=0.$ Then 
\[
\min_{P\in\mathcal{T}_N}\{a_3\}=4(y_N^{(4)})^2-1, \qquad y_N^{(4)}=\sin\frac \pi{N+4}.
\]
The collection of the extremal polynomials $P(z)=z+\sum_{j=2}^N a_j z^j$ forms a one-parameter family, the coefficients of these polynomials are defined by the
formulas 
$$a_{2k+1}=A_{2k+1}(y_N^{(4)}),\qquad a_{2k+2}=\tau A_{2k+2}(y_N^{(4)}),\;
\tau\in[-1,1],\quad k=0,\ldots,(N-2)/2.
$$

\end{theorem}
\bigskip

\begin{remark}
\hspace{1cm}
{}
\\
\hspace*{.4cm}
\textit{(i)} Employing the $y^{(k)}_N, k=1,...,6,$  from Theorem \ref{th:A} the Rogosinski-Szeg\"o estimate \eqref{eq:bounds} can be written as follows:
\begin{align*}
&4(y_N^{(2)})^2-1\le a_3\le4(y_N^{(1)})^2-1, 
\qquad N\textit{ is odd},\ (N+3)/2\textit{ is even},\\
&4(y_N^{(5)})^2-1\le a_3\le4(y_N^{(1)})^2-1, 
\qquad N\textit{ is odd},\ (N+3)/2\textit{ is odd},\\
&4(y_N^{(4)})^2-1\le a_3\le4(y_N^{(3)})^2-1, 
\qquad N\textit{ is even},\ (N+2)/2\textit{ is even},\\
&4(y_N^{(6)})^2-1\le a_3\le4(y_N^{(3)})^2-1, 
\qquad N\textit{ is even},\ (N+2)/2\textit{ is odd},\\
\end{align*}

\textit{(ii)}
For odd $N$, the extremal polynomial $P_N(z)$ is unique and contains only odd powers; in contrast, for even $N$, the extremizer is a sum of a unique polynomial containing
only odd powers and a unique polynomial with even powers multiplied by the 
parameter $\tau\in[-1,1]$. Moreover, if $P_N(z)$ is the unique polynomial for the
case of odd $N$, then the extremal polynomial for the even case can be written as
\begin{equation}\label{eq:oddeven}
P_{N+1}(z)=P_N(z)+\frac12\tau\left(P_N(z)(z+z^{-1})-1\right), 
\end{equation}
$\tau\in[-1,1]$, and 
$\text{extr}_{P\in\mathcal{T}_{N+1}}\{a_3\}
=\text{extr}_{P\in\mathcal{T}_{N}}\{a_3\}$. 

Note that the extremal polynomials turn out to be typically real due to Corollary \ref{cor:2} in Section \ref{sec:5}.
\\

\textit{(iii)} The extremal polynomial for $\max_{P\in \mathcal{T}_4}\{a_3\}$ yields a counterexample to a theorem due to W.C.~Royster and T.J.~Suffridge \cite{RS70}. By Theorem A, d) $y_4^{(3)}=\cos{\pi/4},$ $\max_{P\in \mathcal{T}_4}\{a_3\}=1=a_3=a_1.$ By elementary calculations it follows that $a_2=\tau,$ $a_4=\tau/2,$ $\tau\in[-1,1].$ Thus the extremizer reads
\begin{equation}\label{eq:RS}
P(z)=z+\tau z^2+ z^3 +\frac\tau2z^4,\qquad \tau\in [-1,1].   
\end{equation}
In \cite{RS70} it is claimed that if the extremizer for some $a_k$ reads 
$P(z)=z+\sum_{j=2}^Na_jz^j$ then the polynomial $ R(x)=1+\sum_{j=2}^Na_j U_{j-1}(x)$ has only real roots, located on the segment $[-1,1].$ Without additional assumptions this is not true. $P(z)$ from \eqref{eq:RS} generates $R(x)=4x^2(1+\tau x)$ whose root $x_0=-1/\tau$ lies outside $[-1,1]$ if $0<\tau<1.$\\
\end{remark}

We give some further simple examples resulting from Theorem \ref{th:A}. \\

\newpage

{\bf Examples.}\\

\underline{$N=3$,  the extremizer for $\max_{P\in\mathcal{T}_3}\{a_3\}:$} Since $(y_3^{(1)})^2=1/2$ it follows by Theorem~\ref{th:A}.\textit{a)}, that $\max_{P\in\mathcal{T}_3}\{a_3\}=1,$ hence
\[
P(z)=z+z^3.
\]

\underline{ $N=3$, 
the extremizer for $\min_{P\in\mathcal{T}_3}\{a_3\}$}: Since  $U'_3(x)=24x^2-4$, we have $(y_3^{(5)})^2=1/6$, thus
$\min_{P\in\mathcal{T}_3}\{a_3\}=-1/3$ and hence by Theorem  \ref{th:A},
\[
P(z)=z-1/3z^3.
\]

\underline{$N=4$, 
the extremizer for $\min_{P\in\mathcal{T}_4}\{a_3\}$}: By the above Remark (ii) we have $\min_{P\in\mathcal{T}_4}\{a_3\}=\min_{P\in\mathcal{T}_3}=-1/3$ and
\[
P(z)=z-1/3z^3+\tau/3\cdot(z^2-1/2z^4),\ \tau\in[-1,1].
\]
For larger $N$ the derivation of the extremizers by means of Theorem  \ref{th:A} does not seem appropriate. Therefore, we give a compact form of the extremizers in the shape of pole-free rational functions.

\subsection{Compact representations}

Define rational functions in $z$ by
\[
\hspace{-2cm} G_1(z,y)=32 \frac{y^2(1-y^2)z^5}{(1-z^2)(1+z^4+2(1-2y^2)z^2)^2}
\left\{
\begin{array}{l}
\frac{1-z^{N+5}}{N+5},\ N\textit{odd},\\
\\
\frac{1-z^{N+4}}{N+4},\ N\textit{even},
\end{array}
\right.
\]
\[
\hspace{-7cm} G_2(z,y)=\frac{z+z^3}{1+z^4+2(1-2y^2)z^2},
\]
\[
\hspace{-6.35 cm} G_3(z,y)=\frac{-32y^2(1-y^2)z^7}{(1+z^4+2(1-2y^2)z^2)^2}\times
\]
\[
\hspace{2.8cm} \left\{
\begin{array}{l}
\! \frac{(N+7)^2(1-z^{N+3})+(N+3)^2(1-z^{N+7})-2(N+3)(N+7)(z^2-z^{N+5})}
{(N+3)(N+5)(N+7)(1-z^2)^3 }, N\textit{ odd}, \\ \\
\! \frac{(N+6)^2(1-z^{N+2})+(N+2)^2(1-z^{N+6})-2(N+2)(N+6)(z^2-z^{N+4})}
{(N+2)(N+4)(N+6)(1-z^2)^3}, N\textit{ even},
\end{array}
\right.
\]
\[
\hspace{-3.7 cm} G_4(z,y)=\frac{(z+z^3)(1+z^8-\gamma_1(y)(z^2+z^6)+\gamma_2(y)z^4)}
{(1-z^2)^2(1+z^4+2(1-2y^2)z^2)^2},
\]
where 
\[
\gamma_1(y)=4y^2,
\qquad
\gamma_2(y)=\left\{
\begin{array}{l}
2\left(-\frac{16}{N+5}y^4+4\left(1+\frac{4}{N+5}\right)y^2-1\right),
\ N\textit{ is odd},\\
2\left(-\frac{16}{N+4}y^4+4\left(1+\frac{4}{N+4}\right)y^2-1\right),
\ N\textit{ is even}.\\
\end{array}
\right.
\]

\begin{theorem}\label{th:B}
\hspace{1cm}
{}
\\
\hspace*{.4cm} \textit{a)} Cases $A_\text{max}$, $B_\text{max}$: Let \textit{$N$ be odd},
$y_N^{(1)}$ be the maximum root of the equation $U_{(N+3)/2}(x)=0$. Then the extremizer can be represented by 
\[
P(z)=G_1(z,y_N^{(1)})+G_2(z,y_N^{(1)}). 
\]

\textit{b)} Cases $C_\text{max}$, $D_\text{max}$: Let \textit{$N$ be even} and
$y_N^{(3)}$ be the maximum root of the equation $U_{(N+2)/2}(x)=0$. Then the 
extremizer can be written in the form
\[
P(z)=\left(G_1(z,y_N^{(3)})+G_2(z,y_N^{(3)})\right)\left(1+\tau\frac{z+z^{-1}}2\right)
-\frac{\tau}2,\quad \tau\in[-1,1].
\]

\textit{c)} Case $A_\text{min}$: Let \textit{$N$ be odd} and \textit{$(N+3)/2$ be odd}, $y_N^{(5)}$ be the minimum positive root of the equation
$U'_{(N+3)/2}(x)=0$. Then the extremizer can be described by
\[
P(z)=G_3(z,y_N^{(5)})+G_4(z,y_N^{(5)}).
\]

\textit{d)} Case $B_\text{min}$: Let \textit{$N$ be odd}, \textit{$(N+3)/2$ be even}  and $y_N^{(2)}$ be the minimum positive root of the equation
$U_{(N+3)/2}(x)=0$. Then the extremizer can be represented by
\[
P(z)=G_1(z,y_N^{(2)})+G_2(z,y_N^{(2)}).
\]

\textit{e)} Case $C_\text{min}$: Let \textit{$N$ be even} and \textit{$(N+2)/2$ be odd}, $y_N^{(6)}$ be the minimum positive root of the equation
$U'_{(N+2)/2}(x)=0$. Then the extremizer reads
\[
P(z)=(G_3(z,y_N^{(6)})+G_4(z,y_N^{(6)}))\left(1+\tau\frac{z+z^{-1}}2\right)
-\frac{\tau}2,\quad \tau\in[-1,1].
\]

\textit{f)} Case $D_\text{min}$: Let \textit{$N$ be even}, \textit{$(N+2)/2$ be even}  and $y_N^{(4)}$ be the minimum positive root of the equation 
$U_{(N+2)/2}(x)=0$. Then the extremizer can be described by
\[
P(z)=(G_1(z,y_N^{(4)})+G_2(z,y_N^{(4)}))\left(1+\tau\frac{z+z^{-1}}2\right)
-\frac{\tau}2,\quad \tau\in[-1,1].
\]

\end{theorem}
\medskip

Below are further examples resulting from Theorem \ref{th:B}. \\ \\

\underline{ \textit{(i)} {$N=5$, problem $\max_{P\in\mathcal{T}_5}\{a_3\}$}, i.e. case \textit{a)}}. Since $y_5^{(1)}=\cos\pi/5$, 
$(y_5^{(1)})^2=(1+\cos2\pi/5)/2=(3+\sqrt{5})/8$ it follows that
\[
\max_{P\in\mathcal{T_5}}\{a_3\}=4(y_5^{(1)})^2-1
=1+2\cos\frac{2\pi}{5}=\frac{1+\sqrt{5}}2.
\]
By Theorem~\ref{th:B}.\textit{a)}, the extremizer is given by  
$P(z)=G_1(z,y_5^{(1)})+G_2(z,y_5^{(1)})$, where
\[
G_1(z,y_5^{(1)})=\frac{(5+\sqrt{5})z^5}{10(1-z^2)} 
\frac{1-z^{10}}{\left(1+z^4+\frac{1-\sqrt{5}}2 z^2\right)^2},\qquad
G_2(z,y_5^{(1)})=\frac{z+z^3}{1+z^4+\frac{1-\sqrt{5}}2 z^2}.
\]
After transformations we obtain
\[
P(z)=z+\frac{1+\sqrt{5}}2 z^3+\frac{5+\sqrt{5}}{10}z^5.
\]
\underline{\textit{(ii)} {$N=5$, problem $\min_{P\in\mathcal{T}_5}\{a_3\}$}, i.e. 
case \textit{d)}}. Since $y_5^{(2)}=\sin\pi/10$, 
$(y_5^{(2)})^2=(1-\cos2\pi/5)/2=(3-\sqrt{5})/8$ it follows that
\[
\min_{P\in\mathcal{T}_5}\{a_3\}=4(y_5^{(2)})^2-1=1-2\cos\frac{\pi}5
=\frac{1-\sqrt{5}}{2}.
\]
By Theorem~\ref{th:B}.\textit{d)}, 
$P(z)=G_1(z,y_5^{(2)})+G_2(z,y_5^{(2)})$, where
\[
G_1(z,y_5^{(2)})=\frac{(5-\sqrt{5})z^5}{10(1-z^2)}
\frac{1-z^{10}}{\left(1+z^4+\frac{1+\sqrt{5}}2 z^2\right)^2}\;,\qquad
G_2(z,y_5^{(2)})=\frac{z+z^3}{1+z^4+\frac{1+\sqrt{5}}2 z^2}\,.
\]
After transformations we get 
\[
P(z)=z+\frac{1-\sqrt{5}}2 z^3+\frac{5-\sqrt{5}}{10} z^5.
\]
\underline{ \textit{(iii)} {$N=6$, problem $\max_{P\in\mathcal{T}_6}\{a_3\}$}, 
 i.e. case \textit{b)}}. Since $y_6^{(3)}=y_5^{(1)}=\cos\pi/5$ then, by Remark \textit{(ii)},
$\max_{P\in\mathcal{T}_6}\{a_3\}=(1+\sqrt{5})/2$ and for each $\tau\in[-1,1]$ 
\begin{gather*}
P(z)=z+\frac{1+\sqrt{5}}{2}z^3+\frac{5+\sqrt{5}}{10}z^5
+\tau\left(\frac{3+\sqrt{5}}2 z^2+\frac{5+3\sqrt{5}}5 z^4
+\frac{5+\sqrt{5}}{10}z^6\right)
\end{gather*}
attains the maximum. 

\underline{\textit{(iv)} {$N=6$, problem $\min_{P\in\mathcal{T}_6}\{a_3\}$}, i.e. 
case \textit{f)}}. Analogously, since $y_6^{(4)}=y_5^{(2)}=\sin\pi/10$, we obtain
$\min_{P\in\mathcal{T}_6}\{a_3\}=(1-\sqrt{5})/2$ and for each $\tau\in[-1,1]$ 
\begin{gather*}
P(z)=z+\frac{1-\sqrt{5}}{2}z^3+\frac{5-\sqrt{5}}{10}z^5
+\tau\left(\frac{3-\sqrt{5}}2 z^2+\frac{5-3\sqrt{5}}5 z^4
+\frac{5-\sqrt{5}}{10}z^6\right),\\
\end{gather*}
attains the minimum.

\medskip
Sections \ref{sec:3}, \ref{sec:4}, \ref{sec:6}  of the paper are devoted to the proof of Theorems~\ref{th:A} and
\ref{th:B}. In Section \ref{sec:5} the property "typically real" is proved.

\section{Proof of Theorem~\ref{th:A}}\label{sec:3}
\subsection{A brief survey of the proof}{\it$\mbox{}$ }

\vspace{.2cm}

1. Since the imaginary part of a typically real polynomial on the unit circle is a non-negative
sine-polynomial on $[0,\pi]$, we can reduce the problem to a trigonometric
one. \\

2. By factoring out $\sin t$ the problem is reduced to non-negative cosine polynomials, their coefficients are related to the original ones in \eqref{eq:a-gamma}. \\ 

3. The core of the proof is the application of the 
Fej\'er-Riesz representation to a non-negative cosine-polynomial, which reduces
the problem to the optimization of the positive definite quadratic forms whose 
coefficients are related to the coefficients of the cosine-polynomial by 
formula~\eqref{eq:gamma-delta}.\\

4. The max/min problem for quadratic forms ends up with finding the maximum/minimum eigenvalues of the corresponding matrix pencil. This leads to finding the roots of the determinant of a five-band Toeplitz matrix\\.

5. Its solution also shows that the upper/lower bound of the Rogosinski--Szeg\H{o} estimates is
expressed in terms of the maximum/minimum eigenvalues of the matrix pencil with a subsequent application of the Rayleigh method.
Moreover, we obtain the existence and uniqueness of the extremal polynomials in
the case of odd $N$ and the existence of a one-parameter family of the extremal
polynomials in the case of even $N$.\\

6. Concerning eigenvectors, we have to consider six cases. Their components involve Chebyshev polynomials and
their derivatives and lead to the coefficients of the cosine polynomial in Step 2.\\

7.  The latter, in turn, determines the coefficients of the extremizer. The easiest way to express these coefficients is through recurrence formulas. In the end, these are transformed into formulas giving explicit representations which caused significant technical difficulties.

\subsection{ Fej\'er's idea to convert the extremal problem for $a_3$ to one on eigenvalues and eigenvectors of a matrix pencil}
Let $P(z)=z+\sum_{j=2}^N a_j z^j\in \mathcal{T}_N$. Use the equality
\[
\text{Im}\{P(e^{it})\}=\sin{t}\cdot\tilde{P}(t),\qquad 
\tilde{P}(t)=\gamma_1+2\sum_{j=2}^N \gamma_j\cos(j-1)t.
\]
Here $\tilde{P}(t)$ is a non-negative cosine-polynomial on $[0,\pi]$; the
coefficients $a_1,\ldots,a_N$ and $\gamma_1,\ldots,\gamma_N$ are related by the
bijective relation

\begin{equation}\label{eq:a-gamma}
a_s=\gamma_s-\gamma_{s+2},\quad s=1,\ldots,N.
\end{equation}
On account of \eqref{eq:sine} we have  $a_1 = 1$ and, for conenience, put  
$\gamma_{N+1}=\gamma_{N+2}=0$. It follows from~\eqref{eq:a-gamma} that 
$\gamma_1-\gamma_3=1$, $a_3=\gamma_3-\gamma_5$. 

The polynomial $\tilde{P}(t)$ is non-negative, and by the Fej\'er--Riesz theorem
it can be represented in the form 
$\tilde{P}(t)=|\delta_1+\delta_2 e^{it}+\ldots+\delta_N e^{i(N-1)t}|^2$, whence
\begin{equation}\label{eq:gamma-delta}
\gamma_s=\sum_{j=1}^{N-s+1} \delta_j \delta_{j+s-1},\quad s=1,\ldots,N.
\end{equation}
Then 
\[
a_3=\gamma_3-\gamma_5=\sum_{j=1}^{N-2} \delta_j \delta_{j+2}-\sum_{j=1}^{N-4} \delta_j \delta_{j+4},
\quad
1=\gamma_1-\gamma_3=\sum_{j=1}^N \delta_j^2-\sum_{j=1}^{N-2} \delta_j \delta_{j+2}.
\] 

The quadratic form $\sum_{j=1}^N \delta_j^2-\sum_{j=1}^{N-2} \delta_j
\delta_{j+2}$ is positive definite (see~\cite[Lemma~A.1]{DST24}),
therefore~\cite{MW68}, 
\begin{align*}
\min\left\{\sum_{j=1}^{N-2}\delta_j\right.&\delta_{j+2}\left.
-\sum_{j=1}^{N-4}\delta_j\delta_{j+4}:
\sum_{j=1}^N \delta_j^2-\sum_{j=1}^{N-2}\delta_j\delta_{j+2}=1\right\} \\
&\le a_3 \le
\max\left\{\sum_{j=1}^{N-2}\delta_j\delta_{j+2}
-\sum_{j=1}^{N-4}\delta_j\delta_{j+4}:
\sum_{j=1}^N \delta_j^2-\sum_{j=1}^{N-2}\delta_j\delta_{j+2}=1\right\}.
\end{align*}

Define two symmetric matrices $A$ and $B$ of order $N\times N$ associated to the
quadratic forms 
\[
\sum_{j=1}^{N-2}\delta_j\delta_{j+2}-\sum_{j=1}^{N-4}\delta_j\delta_{j+4}
\quad\text{and}\quad 
\sum_{j=1}^N \delta_j^2-\sum_{j=1}^{N-2}\delta_j\delta_{j+2},
\]
\[
A=
\begin{pmatrix}
0 & 0 & 1/2 & 0 & -1/2 & \ldots \\
0 & 0 & 0 & 1/2 & 0 & \ldots \\
1/2 & 0 & 0 & 0 & 1/2 & \ldots \\
0 & 1/2 & 0 & 0 & 0 & \ldots \\
-1/2 & 0 & 1/2  & 0 & 0 & \ldots \\
\ldots & \ldots & \ldots & \ldots & \ldots & \ldots
\end{pmatrix},
\;
B=
\begin{pmatrix}
1 & 0 & -1/2 & 0 & \ldots \\
0 & 1 & 0 & -1/2 & \ldots \\
-1/2 & 0 & 1 & 0 & \ldots \\
0 & -1/2 & 0 & 1 & \ldots \\
0 & 0 & -1/2  & 1 & \ldots \\
\ldots & \ldots & \ldots & \ldots & \ldots
\end{pmatrix}
\]
respectively.

Note that the matrix $A$ has two non-zero diagonals above the main diagonal: on the 
third diagonal all elements are equal to $1/2$, on the fifth one they all equal
$-1/2$. Similarly, matrix $B$ has ``units'' on the main diagonal and $-1/2$
on the third diagonal. Below the main diagonal, the elements are arranged 
symmetrically.

Let $\lambda_1\le\ldots\le\lambda_N$ be the roots of the equation 
\begin{equation}\label{eq:char-eq}
\det(A-\lambda B)=0.
\end{equation}

The numbers $\lambda_1,\ldots,\lambda_N$ are the eigenvalues of the matrix pencil $\{A-\lambda B,\lambda\in\mathbb{C}\}$. Then 
\[
\lambda_1\le a_3\le\lambda_N.
\] 

To find the extremizers, we need to know the eigenvectors corresponding to the
eigenvalues $\lambda_1$ and $\lambda_N$, that is, we have to find non-trivial solutions of the
equations $(A-\lambda_N B)Z=0$ and $(A-\lambda_1 B)Z=0$. 

Let the vector $Z^{(0)}=(z_1^{(0)},\ldots,z_N^{(0)})^T$ be a solution to the first
equation (the sign $T$ denotes transposition). Then by 
formulas~\eqref{eq:gamma-delta} and \eqref{eq:a-gamma}, the coefficients of the
extremizer for the problem $\max\{a_3\}$ are defined by

\begin{equation}\label{eq:gamma-0}
\gamma_s^{(0)}=\sum_{j=1}^{N-s+1} z_j^{(0)} z_{j+s-1}^{(0)},\quad
s=1,\ldots,N,
\end{equation}

\begin{equation}\label{eq:a-0}
a_s^{(0)}=\frac{\gamma_s^{(0)}-\gamma_{s+2}^{(0)}}
{\gamma_1^{(0)}-\gamma_3^{(0)}},\quad
s=1,\ldots,N,
\end{equation}
where we assume that $\gamma_{N+1}^{(0)}=\gamma_{N+2}^{(0)}=0$. Because $Z^{(0)}$ is a solution of a homogeneous system, we have to normalize $a_s^{(0)}$, as one can see by comparing \eqref{eq:a-gamma} and \eqref{eq:a-0}.
The coefficients of the extremizer for the problem $\min\{a_3\}$ are determined
analogously. Thus, the problem has been reduced to the determination of the eigenvalues of the matrix pencil
$\{A-\lambda B,\lambda\in\mathbb{C}\}$ and of their eigenvectors.
\subsection{Computing the determinant $\det(A-\lambda B)$}
Set $\lambda=4x^2-1$ and consider the $N\times N$ matrix
\begin{align*}
&\Phi_N (x)=A-\lambda B \\
&\quad=\begin{pmatrix}
-4x^2+1 & 0 & 2x^2 & 0 & -1/2 & 0 & \ldots \\
0 & -4x^2+1 & 0 & 2x^2 & 0 & -1/2 & \ldots \\
2x^2 & 0 & -4x^2+1 & 0 & 2x^2 & 0 & \ldots \\
0 & 2x^2 & 0 & -4x^2+1 & 0 & 2x^2 & \ldots \\
-1/2 & 0 & 2x^2 & 0 & -4x^2+1 & 0 & \ldots \\
0 & -1/2 & 0 & 2x^2 & 0 & -4x^2+1 & \ldots \\
\ldots & \ldots & \ldots & \ldots & \ldots & \ldots & \ldots  
\end{pmatrix}.
\end{align*}

\begin{ntheorem}\label{th:Delta_N}
If we set $\Delta_N=\det(A-\lambda B),$ $\lambda=4x^2-1,$ then
\begin{equation}\label{eq:det}
\Delta_N=\left\{
\begin{array}{l}\displaystyle
-\frac{1}{2^{N+4}x^2}U_\frac{N+1}{2}(x)U'_\frac{N+1}{2}(x)U_\frac{N+3}{2}
U'_\frac{N+3}{2}(x),\quad \textit{ $N$ odd},\\
\displaystyle
\frac{1}{2^{N+4}x^2}\left[U_\frac{N+2}{2}(x)U'_\frac{N+2}{2}(x)\right]^2,\quad \textit{ $N$ even}.
\end{array}
\right.
\end{equation}
\end{ntheorem}
\bigskip


\begin{proof}
Let us perform the following transformations on the columns and rows of the 
matrix $\Phi_N(x)$. First, we swap the third and the second column, the third and
the second row. Then, swap consecutively the fifth and the fourth column, the 
fourth and the third column; the same with the rows. We use this scheme 
$(N-1)/2$ times. The last chain of column and row transformations is 
$N\to N-1\to\ldots\to(N+1)/2$ (when \textit{$N$ is odd}), 
$N-1\to N-2\to\ldots\to N/2$ (when \textit{$N$ is even}). Each of the
transformation chains does not change the determinant of the matrix. As a result,
we get that the matrix $\Phi_N(x)$ will be transformed into the matrix
\[
\begin{pmatrix}
D_\frac{N+1}{2} & O \\
O & D_\frac{N-1}{2}
\end{pmatrix}
\ (N\textit{ is odd})
\quad\text{or}\quad
\begin{pmatrix}
D_\frac{N}{2} & O \\
O & D_\frac{N}{2}
\end{pmatrix}
\ (N\textit{ is even}),
\]
where $O$ is the zero matrix and the $N\times N$ matrix $D_N$ is defined
\begin{equation}\label{DN}
D_N (x)=\begin{pmatrix}
-4x^2+1 & 2x^2 & -1/2 & 0 & \ldots \\
2x^2 & -4x^2+1 & 2x^2 & -1/2 & \ldots \\
-1/2& 2x^2 & -4x^2+1 & 2x^2 & \ldots \\
0 & -1/2 & 2x^2 & -4x^2+1 & \ldots \\
\ldots & \ldots & \ldots & \ldots & \ldots 
\end{pmatrix}
\end{equation}
 Now invoke the following result from
\cite{DSS22,DSS23}
\begin{equation}\notag
\det D_N(x)=\frac{(-1)^N}{2^{N+2}x}U_{N+1}(x)U'_{N+1}(x).
\end{equation}

Thus, by applying the rule for calculating the determinant of a block-diagonal matrix, we arrive at the
conclusion of the theorem.
\end{proof}

\subsection{Extremal values of the coefficient $a_3$}
Using formulas~\eqref{eq:det}, we easily obtain the extremal values of the 
coefficient $a_3$ in all the cases considered (Lemma~\ref{le:A4}).

\textit{Cases $A_\text{max}$, $B_\text{max}$ ($N$ is odd).} 
$\lambda_N^{(1)}=4(y_N^{(1)})^2-1$, where $y_N^{(1)}$ is the maximum root of the
equation $U_{(N+3)/2}(x)=0$, i.e., $\lambda_N^{(1)}=4\cos^2{2\pi/(N+5)}-1$.

\textit{Cases $C_\text{max}$, $D_\text{max}$ ($N$ is even).} 
$\lambda_N^{(3)}=4(y_N^{(3)})^2-1$, where $y_N^{(3)}$ is the maximum root of the
equation $U_{(N+2)/2}(x)=0$, i.e., $\lambda_N^{(3)}=4\cos^2{2\pi/(N+4)}-1$.

\textit{Case $A_\text{min}$ ($N$ is odd, $(N+3)/2$ is odd).} 
$\lambda_N^{(5)}=4(y_N^{(5)})^2-1$, where $y_N^{(5)}$ is the minimum root of the
equation $U'_{(N+3)/2}(x)=0$.

\textit{Case $B_\text{min}$ ($N$ is odd, $(N+3)/2$ is even).} 
$\lambda_N^{(2)}=4(y_N^{(2)})^2-1$, where $y_N^{(2)}$ is the minimum root of the
equation $U_{(N+3)/2}(x)=0$, i.e., $\lambda_N^{(2)}=4\sin^2{\pi/(N+5)}-1$.

\textit{Case $C_\text{min}$ ($N$ is even, $(N+2)/2$ is odd).} 
$\lambda_N^{(6)}=4(y_N^{(6)})^2-1$, where $y_N^{(6)}$ is the minimum root of the
equation $U'_{(N+2)/2}(x)=0$.

\textit{Case $D_\text{min}$ ($N$ is even, $(N+2)/2$ is even).} 
$\lambda_N^{(4)}=4(y_N^{(4)})^2-1$, where $y_N^{(4)}$ is the minimum root of the
equation $U_{(N+2)/2}(x)=0$, i.e., $\lambda_N^{(4)}=4\sin^2{\pi/(N+4)}-1$.

In the cases $A_\text{max}$, $B_\text{max}$, $B_\text{min}$, $C_\text{max}$, 
$D_\text{max}$, $D_\text{min}$ it is easy to check that the formulas for the 
bounds of the coefficient $a_3$ coincide with the corresponding formulas 
in~\eqref{eq:bounds}.

Consider the case $A_\text{min}$. After the substitution
$x=\cos{t}$ the equation $U'_{(N+3)/2}(x)=0$ turns into the equation (see. \cite{DSS22,DSS23})
\[
\frac{N+7}{2}\sin\frac{N+3}{2}t-\frac{N+3}{2}\sin\frac{N+7}{2}t=0.
\]
Set $n'=(N-1)/2$ and substitute  $\vartheta/2=t+\pi/2$  to transform this equation into 
\[
(n'+4)\cos\frac{(n'+2)\vartheta}{2}+(n'+2)\cos\frac{(n'+4)\vartheta}{2}=0.
\]
Thus, if one needs in the case $A_\text{min}$ the minimum positive root of $U_{(N+3)/2}(x)=0,$  this
coincides with the condition made in the statement of~\eqref{eq:bounds}.

The case $C_\text{min}$ is treated similarly.

\subsection{Finding the eigenvectors of the matrix pencil 
$\{A-\lambda B,\lambda\in\mathbb{C}\}$}
\subsubsection{Cases a), d) ($A_\text{max}$, $B_\text{max}$, $B_\text{min}$)}
\begin{ntheorem}\label{th:2}
Let \textit{$N$ be odd}. The eigenvector $Z$ corresponding to the eigenvalue 
$\lambda_N^{(1)}=4\cos^2 2\pi/(N+5)-1$ of the matrix pencil 
$\{A-\lambda B, \lambda\in\mathbb{C}\}$ is equal to 
$Z=Z^{(0)}(\cos2\pi/(N+5))$, where
\[
Z^{(0)}(x)=\left\{\zeta_1(x),0,\zeta_2(x),0,\ldots,\zeta_\frac{N+1}{2}(x)\right\},
\]
\[
\zeta_k(x)=U_{k-1}(x)U_k(x),\quad k=1,\ldots,\frac{N+1}{2}.
\]
\end{ntheorem}
\begin{proof}
Consider the vector $\hat{Z}^{(0)}(x)=\{U_{k-1}(x)U_k(x)\}_{k=1}^{(N+1)/2}$.
By Lemma~\ref{le:A5}, 
\[
D_\frac{N+1}{2}(\hat{x})\cdot\hat{Z}^{(0)}(\hat{x})=\mathbf{0},
\]
where $\hat{x}=\cos2\pi/(N+5)$, $\mathbf{0}$ is the zero vector, and the matrix
$D_N(x)$ is defined by \eqref{DN}. The components of the vector 
$\Phi_N(x)Z^{(0)}(x)$ with odd indices coincide with the components of the 
vector $D_{(N+1)/2}(x)\cdot\hat{Z}^{(0)}(x)$ (components are compared in 
sequential order), and those with even indices are equal to zero. Hence,
$\Phi_N(\hat{x})Z^{(0)}(\hat{x})=\mathbf{0}$, The theorem is proved.
\end{proof}
\begin{ntheorem}\label{th:3}
Let \textit{$N$ be odd, $(N+3)/2$ be even}. The eigenvector $Z$ corresponding to
the eigenvalue $\lambda_N^{(2)}=4\sin^2\pi/(N+5)-1$ of the matrix pencil 
$\{A-\lambda B, \lambda\in\mathbb{C}\}$ is equal to 
$Z=Z^{(0)}(\sin\pi/(N+5))$, where the vector $Z^{(0)}(x)$ is defined as in the 
previous theorem.
\end{ntheorem}
\begin{proof}
Note that 
\[
\cos\left(\frac{\pi}{N+5}\frac{N+3}{2}\right)=\sin\frac{\pi}{N+5},
\]
that is, $\lambda_N^{(2)}=4\sin^2\frac\pi{N+5}-1$. Then the proof is carried out 
similarly to the previous proof. 
\end{proof}
Thus, in the cases $A_\text{max}$, $B_\text{max}$, $B_\text{min}$, the 
eigenvectors form a one-dimensional subspace. 
\subsubsection{Cases b), f) ($C_\text{max}$, $D_\text{max}$, $D_\text{min}$)}
In these cases, the largest positive eigenvalue of the matrix pencil
$\{A-\lambda B, \lambda\in\mathbb{C}\}$ equals 
$\lambda_N^{(3)}=4\cos^2 2\pi/(N+4)-1$, and the smallest one is 
$\lambda_N^{(4)}=4\sin^2 2\pi/(N+4)-1$, where these eigenvalues have second
multiplicity. Thus to each of them, there will correspond two linearly
independent eigenvectors (clearly, with the corresponding normalization).
\begin{ntheorem}\label{th:4}
Let \textit{$N$ be even}. The linearly independent eigenvectors $Z_1$ and $Z_2$
corresponding to the eigenvalue $\lambda_N^{(3)}=4\cos^2 2\pi/(N+4)-1$ of the 
matrix pencil $\{A-\lambda B,\lambda\in\mathbb{C}\}$ are 
\[
Z_1=Z_1^{(0)}\left(\cos\frac{2\pi}{N+4}\right)
\quad\text{and}\quad
Z_2=Z_2^{(0)}\left(\cos\frac{2\pi}{N+4}\right),
\]
where
\[
Z_1^{(0)}(x)
=\left\{\zeta_1(x),0,\zeta_2(x),0,\ldots,\zeta_\frac{N}{2}(x),0\right\},
\]
\[
Z_2^{(0)}(x)=\left\{0,\zeta_1(x),0,\ldots,\zeta_\frac{N}{2}(x)\right\},
\]
\[
\zeta_k(x)=U_{k-1}(x)U_k(x),\quad k=1,\ldots,\frac{N}{2}.
\]
\end{ntheorem}

\begin{proof}
Consider the vector $\hat{Z}^{(0)}(x)=\{U_{k-1}(x)U_k(x)\}_{k=1}^{N/2}$.
By Lemma~\ref{le:A5}, 
\[
D_\frac{N}{2}(\hat{x})\cdot\hat{Z}^{(0)}(\hat{x})=\mathbf{0},
\]
where $\hat{x}=\cos2\pi/(N+4)$, $\mathbf{0}$ is the zero vector, and the matrix
$D_N(x)$ is defined by \eqref{DN}. The components of the vector 
$\Phi_N(x)Z_1^{(0)}(x)$ with odd indices coincide with the components of the 
vector $D_{N/2}(x)\cdot\hat{Z}^{(0)}(x)$ (components are compared in 
sequential order), and those with even indices are equal to zero. For the 
vector $\Phi_N(x)Z_2^{(0)}(x)$, we have an analogous situation: its components
with even indices coincide with the components of the vector 
$D_{N/2}(x)\cdot\hat{Z}^{(0)}(x)$, while those with odd indices equal zero. Hence,
$\Phi_N(\hat{x})Z_1^{(0)}(\hat{x})=\mathbf{0}$ and 
$\Phi_N(\hat{x})Z_2^{(0)}(\hat{x})=\mathbf{0}$. The theorem is proved.  
\end{proof}
Similarly, we obtain the following theorem.
\begin{ntheorem}
Let \textit{$B$ be even, $(N+2)/2$ be even}. The linearly independent
eigenvectors $Z_1$ and $Z_2$ corresponding to the eigenvalue
$\lambda_N^{(4)}=4\sin^2\pi/(N+4)-1$ of the matrix pencil 
$\{A-\lambda B,\lambda\in\mathbb{C}\}$ equal
\[
Z_1=Z_1^{(0)}\left(\sin\frac{\pi}{N+4}\right)
\quad\text{and}\quad
Z_2=Z_2^{(0)}\left(\sin\frac{\pi}{N+4}\right),
\]
where the vectors $Z_1^{(0)}(x)$ and $Z_2^{(0)}(x)$ are defined as in the 
previous theorem.
\end{ntheorem}
\begin{proof}
It suffices to note that 
\[
\cos\left(\frac{\pi}{N+4}\frac{N+2}{2}\right)=\sin\frac{\pi}{N+4},
\]
that is, $\lambda_N^{(4)}=4\sin^2\frac\pi{N+4}-1$.
\end{proof}
Thus, in the cases $C_\text{max}$, $D_\text{max}$, $D_\text{min}$, the 
eigenvectors form a two-dimensional subspace.
\subsubsection{Case c) ($A_\text{min}$)}
\begin{ntheorem}\label{th:6}
Let \textit{$N$ be odd, $(N+3)/2$ be odd}; $y_N^{(5)}$ be the minimum positive
root of the equation $U'_{(N+3)/2}(x)=0$. The eigenvector $Z$ corresponding to
the eigenvalue $\lambda_N^{(5)}=4y_N^{(5)}-1$ of the matrix pencil 
$\{A-\lambda B, \lambda\in\mathbb{C}\}$ is equal to 
$Z=Z^{(1)}(y_N^{(5)})$, 
\[
Z^{(1)}(x)=\left\{\zeta_1(x),0,\zeta_2(x),0,\ldots,\zeta_\frac{N+1}{2}(x)\right\},
\]
Let the $\zeta_k(x)$ be given for $k=1,\ldots,(N+1)/2$ by
\[
\zeta_k(x)=U_{k-1}(x)U_k(x)-U_{\frac{N+1}{2}-k}(x)U_{\frac{N+3}{2}-k}(x)
+\left(1-\frac{4k}{N+3}\right)U_\frac{N+1}{2}(x)U_\frac{N+3}{2}(x).
\]
\end{ntheorem}
\begin{proof}
Consider the vector $\hat{Z}^{(1)}(x)=\{\zeta_k(x)\}_{k=1}^{(N+1)/2}$.
By Lemma~\ref{le:A6}, 
\[
D_\frac{N+1}{2}(\nu_1)\cdot\hat{Z}^{(1)}(\nu_1)=\mathbf{0},
\]
where $\mathbf{0}$ is the zero vector, $\nu_1=y_N^{(5)}$, and the matrix
$D_N(x)$ is defined by \eqref{DN}. The components of the vector 
$\Phi_N(x)Z^{(1)}(x)$ with odd indices coincide with the components of the 
vector $D_{(N+1)/2}(x)\cdot\hat{Z}^{(1)}(x)$ (components are compared in 
sequential order), and those with even indices are equal to zero. Hence,
$\Phi_N(\nu_1)Z^{(0)}(\nu_1)=\mathbf{0}$, The theorem is proved. 
\end{proof}
\subsubsection{Case e) ($C_\text{min}$)}
\begin{ntheorem}\label{th:7}
Let \textit{$N$ be even, $(N+2)/2$ be odd}; $y_N^{(6)}$ be the minimum positive
root of the equation $U'_{N/2}(x)=0$. The eigenvectors $Z_1$ and $Z_2$
corresponding to the eigenvalue $\lambda_N^{(5)}=4y_N^{(5)}-1$ of the 
matrix pencil $\{A-\lambda B,\lambda\in\mathbb{C}\}$ equal
$Z_1=Z_1^{(1)}(y_N^{(6)})$ and $Z_2=Z_2^{(1)}(y_N^{(6)})$. To describe their components  we set for 
$k=1,\ldots, N/2$
\[
\zeta_k(x)=U_{k-1}(x)U_k(x)-U_{\frac{N}{2}-k}(x)U_{\frac{N+2}{2}-k}(x)
+\left(1-\frac{4k}{N+2}\right)U_\frac{N}{2}(x)U_\frac{N+2}{2}(x),
\]
then

\[
Z_1^{(1)}(x)
=\left\{\zeta_1(x),0,\zeta_2(x),0,\ldots,\zeta_\frac{N}{2}(x),0\right\},
\]
\[
Z_2^{(1)}(x)=\left\{0,\zeta_1(x),0,\zeta_2(x),0,\ldots,
\zeta_\frac{N}{2}(x)\right\}.
\]
\end{ntheorem}
\begin{proof}[Proof\nopunct]
is analogous to Theorems~\ref{th:4} and \ref{th:6}.
\end{proof}
\subsection{Recurrence formulas for the coefficients of the extremal 
polynomials.}
In the case when \textit{$N$ is odd}, the eigenvectors form a one-dimensional 
subspace and will determine the coefficients of the corresponding extremizers 
in a unique way. In the case when \textit{$N$ is even}, the eigenvectors form a
two-dimensional subspace, therefore, the coefficients of the extremizers will
form a one-parameter set, which, as we show below, will turn out to be a segment.
\subsubsection{$N$ is odd} 
Let us use formulas~\eqref{eq:a-gamma} and \eqref{eq:gamma-delta}:

\begin{equation}\label{eq:gamma-z}
\gamma_s=\sum_{j=1}^{N-s+1} z_j z_{j+s-1},\quad
s=1,\ldots,N,
\end{equation}

\begin{equation}\label{eq:a-frac.gamma}
a_s=\frac{\gamma_s-\gamma_{s+2}}
{\gamma_1-\gamma_3},\quad
s=1,\ldots,N.
\end{equation}

Here the vector $\{z_1,\ldots,z_N\}$ is one of the eigenvectors defined in
Theorems~\ref{th:2}, \ref{th:3}, \ref{th:6}.

\textit{Case a) ($A_\text{max}$, $B_\text{max}$).} By 
formulas~\cite[pp.~20, 21]{DSS22} and \cite[pp.~7, 8]{DSS23}, we derive the 
relations
\[
\gamma_{2j-1}-\gamma_{2j+3}=y_N^{(1)}U_{j-1}(y_N^{(1)}) 
U'_{\frac{N+1}{2}-j+1}(y_N^{(1)}),
\]
\[
\gamma_{2j}-\gamma_{2j+4}=0,\quad
j=1,\ldots,\frac{N+1}{2},
\]
where we recall $y_N^{(1)}=\cos2\pi/(N+5)$ (we assume that 
$\gamma_{N+1}=\gamma_{N+2}=\gamma_{N+3}=\gamma_{N+4}=0$). Note that
\[
\gamma_1-\gamma_5=y_N^{(1)} U'_\frac{N+1}{2}(y_N^{(1)})=\frac{N+5}{2}
\frac{(y_N^{(1)})^2}{1-(y_N^{(1)})^2}.
\]
We need to calculate $\gamma_s-\gamma_{s+2}$, $s=1,\ldots,N-2$.

For $s=1$,
\[
a_3=4(y_N^{(1)})^2-1=\frac{\gamma_3-\gamma_5}{\gamma_1-\gamma_3}
=\frac{\gamma_3-\gamma_5+\gamma_1-\gamma_1}{\gamma_1-\gamma_3}
=-1+\frac{\gamma_1-\gamma_5}{\gamma_1-\gamma_3},
\]
whence
\[
\gamma_1-\gamma_3=\frac{N+5}{8}\frac1{1-(y_N^{(1)})^2}.
\]
Further,
\[
a_{2k-1}=\frac{\gamma_{2k-1}-\gamma_{2k+1}+\gamma_{2k+3}-\gamma_{2k+3}}
{\gamma_1-\gamma_3}=-a_{2k+1}+\frac{\gamma_{2k-1}-\gamma_{2k+3}}
{\gamma_1-\gamma_3}.
\]
Then for $k=1,\ldots, (N-1)/2$ there holds

\begin{equation}\label{a:case-a}
a_{2k+1}=-a_{2k-1}+\frac{8(1-(y_N^{(1)})^2)y}{N+5} U_{k-1}(y_N^{(1)})
U'_{\frac{N+3}{2}-k}(y_N^{(1)}),\quad 
a_1=1, a_{2k}=0.
\end{equation}
\textit{Case d) ($B_\text{min}$).}
This case is treated similarly. Let $y_N^{(2)}=\sin\pi/(N+5)$. Then
\[
\gamma_{2j-1}-\gamma_{2j+3}=(-1)^{\frac{N-1}{4}}y_N^{(2)}U_{j-1}(y_N^{(2)}) 
U'_{\frac{N+1}{2}-j+1}(y_N^{(2)}),
\]
\[
\gamma_{2j}-\gamma_{2j+4}=0,\quad
j=1,\ldots, (N+1)/2,
\]
\[
\gamma_1-\gamma_5=(-1)^{\frac{N-1}{4}}y_N^{(2)} U'_\frac{N+1}{2}(y_N^{(2)})=\frac{N+5}{2}
\frac{(y_N^{(2)})^2}{1-(y_N^{(2)})^2}.
\]
(here again $\gamma_{N+1}=\gamma_{N+2}=\gamma_{N+3}=\gamma_{N+4}=0$). 

Further, $a_1=1$ and for $k=1,\ldots, (N-1)/2$ we have $a_{2k}=0,$ thus
\begin{equation}\label{a:case-b}
\displaystyle
a_{2k+1}=-a_{2k-1}+(-1)^{\frac{N-1}{4}}
\frac{8(1-(y_N^{(2)})^2)y}{N+5} U_{k-1}(y_N^{(2)})
U'_{\frac{N+3}{2}-k}(y_N^{(2)}).
\end{equation}

\textit{Case c) ($A_\text{min}$).} 
Let $y_N^{(5)}$ be the minimum positive root of $U'_{(N+3)/2}(x)=0.$
 Choose $z_k=\zeta_K(y_N^{(5)})$, $\zeta_k(x)$ being defined as in Theorem~\ref{th:6}. Then
by formulas~\eqref{eq:gamma-z} and \eqref{eq:a-frac.gamma} we have
\begin{equation}\label{a:case-e}
a_{2k-1}=\frac{\sum\limits_{j=1}^{(N+3)/2-k}z_j z_{j+k-1}
-\sum\limits_{j=1}^{(N+1)/2-k}z_j z_{j+k}}
{\sum\limits_{j=1}^{(N+1)/2}z_j z_j
-\sum\limits_{j=1}^{(N-1)/2}z_j z_{j+1}}, \quad
k=1,\ldots,\frac{N+1}{2}.
\end{equation}

Thus, in the case of odd $N$, the extremizers are determined uniquely,
the coefficients with even indices being zero; the coefficients with odd indices for 
the problems $A_\text{max}$, $B_\text{max}$, $B_\text{min}$, $A_\text{min}$ are 
defined by formulas~\eqref{a:case-a}, \eqref{a:case-b} and \eqref{a:case-e}, 
respectively. 
\subsubsection{$N$ is even} 
When $N$ is even, the coefficients of the odd $z$-powers are determined similarly to the
case of odd $N$. However, the coefficients associated to even powers will not necessarily be
zero; besides, they will linearly depend on a parameter. 

\textit{Case b) ($C_\text{max}$, $D_\text{max}$).}
Let $y_N^{(3)}=\cos2\pi/(N+4)$, $Z_1=Z_1^{(0)}(y_N^{(3)})$, 
$Z_2=Z_2^{(0)}(y_N^{(3)})$, where $Z_1^{(0)}(x)$
and $Z_2^{(0)}(x)$ are defined as in Theorem~\ref{th:4}. The values $\gamma_s$
$(s=1,\ldots,N)$ are determined by formulas~\eqref{eq:gamma-z}, where
\[
z_{2k-1}=\alpha\zeta_k(y_N^{(3)}), \qquad z_{2k}=\beta\zeta_k(y_N^{(3)}), \quad k=1,\ldots,N/2,
\]
and $\alpha$, $\beta$ are arbitrary constants. Then it follows 
from~\cite[pp.~20, 21]{DSS22} and \cite[pp.~7, 8]{DSS23} that 
\[
\hspace{-.4cm} \gamma_{2j-1}-\gamma_{2j+3}=(\alpha^2+\beta^2)y_N^{(3)}U_{j-1}(y_N^{(3)})
U'_{\frac{N+2}{2}-j}(y_N^{(3)}),\qquad j=1,\ldots, N/2,
\]
\[
\gamma_{2j}-\gamma_{2j+4}=\alpha\beta y_N^{(3)}
\left(U_{j-1}(y_N^{(3)})U'_{\frac{N+2}{2}-j}(y_N^{(3)})+U_j(y_N^{(3)})
U'_{\frac{N}{2}-j}(y_N^{(3)})\right),
\]
\[
\gamma_1-\gamma_3=(\alpha^2+\beta^2)\frac{N+4}{8}\frac1{1-(y_N^{(3)})^2},
\]
(we assume that $\gamma_{N+1}=\gamma_{N+2}=\gamma_{N+3}=\gamma_{N+4}=0$).

The coefficients with odd indices are determined uniquely: $a_1=1$ and for $k=1,\ldots, (N-2)/2$ by
\begin{equation}\label{a-odd:case-c}
a_{2k+1}=-a_{2k-1}+\frac{8(1-(y_N^{(3)})^2)y}{N+4} U_{k-1}(y_N^{(3)})
U'_{\frac{N+2}{2}-k}(y_N^{(3)}). \\
\end{equation}
We have to find the coefficients with even indices: 
\[
a_{2k}=\frac{\gamma_{2k}-\gamma_{2k+2}+\gamma_{2k+4}-\gamma_{2k+4}}
{\gamma_1-\gamma_3}=-a_{2k+2}+\frac{\gamma_{2k}-\gamma_{2k+4}}
{\gamma_1-\gamma_3}.
\]
Define $\tau=2\alpha\beta/(\alpha^2+\beta^2)\in[-1,1]$. Then $a_0=0$ and for $k=1,\ldots,(N-2)/2$ there holds  
\begin{equation}\label{a-even:case-c}
\begin{array}{l}
\begin{aligned}
a_{2k+2}=-a_{2k}&+\tau\frac{4(1-(y_N^{(3)})^2)y_N^{(3)}}{N+4} \\
&\times\left(U_{k-1}(y_N^{(3)}) U'_{\frac{N+2}{2}-k}(y_N^{(3)})
+U_k(y)U'_{\frac{N}2-k}(y_N^{(3)})\right), 
\end{aligned} \\
\end{array}
\end{equation}
Note that \eqref{a-even:case-c} easily implies that 
$a_2=2\tau(y_N^{(3)})^2=2\tau\cos^2 2\pi/(N+4)$. 

\textit{Case f) ($D_\text{min}$).}
This case is treated similarly. Recall $y_N^{(4)}=\sin\pi/(N+4)$. Then $a_0=0,\,a_1=1,$ and for  $k=1,\ldots,(N-2)/2$ there holds
\begin{equation}\label{a:case-d}
\begin{array}{l}
\begin{aligned}
a_{2k+2}=-a_{2k}&+(-1)^{\frac{N-1}{4}}\tau
\frac{4(1-(y_N^{(4)})^2)y_N^{(4)}}{N+4} \\
&\times\left(U_{k-1}(y_N^{(4)}) U'_{\frac{N+2}{2}-k}(y_N^{(4)})
+U_k(y_N^{(4)})U'_{\frac{N}2-k}(y_N^{(4)})\right), 
\end{aligned}\\
\displaystyle
a_{2k+1}=-a_{2k-1}+(-1)^{\frac{N-2}{4}}\frac{8(1-(y_N^{(4)})^2)y}{N+4} 
U_{k-1}(y_N^{(4)}) U'_{\frac{N+2}{2}-k}(y_N^{(4)}). \\
\end{array}
\end{equation}

\textit{Case e) ($C_\text{min}$).}
Let $y_N^{(6)}$ be the minimum positive root of 
$U'_{(N+2)/2}(x)=0$,
\[
z_{2k-1}=\alpha\zeta_k(y_N^{(6)}), \qquad
z_{2k}=\beta\zeta_k(y_N^{(6)}), \qquad
k=1,\ldots,N/2,
\]
with $\zeta_k(x)$ defined as in Theorem~\ref{th:7}. Then $\gamma_N=\alpha\beta z_1 z_\frac{N}{2},$ and for  $k=1,\ldots,(N-2)/2$ one has
\[
\begin{array}{l}
\displaystyle
\gamma_{2k-1}=(\alpha^2+\beta^2)\sum_{j=1}^{(N+2)/2-k}z_j z_{j+k-1},\\
\displaystyle
\gamma_{2k}=\alpha\beta\left(\sum_{j=1}^{(N+2)/2-k}z_j z_{j+k-1}
+\sum_{j=1}^{N/2-k}z_j z_{j+k}\right).\\
\end{array}
\]
By formulas~\eqref{eq:gamma-z}, \eqref{eq:a-frac.gamma} we get for $\tau\in[-1,1],\,k=1,\ldots, N/2$
\begin{equation}\label{a:case-f}
\left\{
\begin{array}{l}
\displaystyle
a_{2k-1}=\frac{\sum\limits_{j=1}^{(N+2)/2-k}z_j z_{j+k-1}
-\sum\limits_{j=1}^{N/2-k}z_j z_{j+k}}
{\sum\limits_{j=1}^{N/2}z_j z_j
-\sum\limits_{j=1}^{(N-2)/2}z_j z_{j+1}}, \\
\displaystyle
a_{2k}=\tau\frac{\sum\limits_{j=1}^{(N+2)/2-k}z_j z_{j+k-1}
-\sum\limits_{j=1}^{(N-2)/2-k}z_j z_{j+k+1}}
{\sum\limits_{j=1}^{N/2}z_j z_j
-\sum\limits_{j=1}^{(N-2)/2}z_j z_{j+1}}. \\
\end{array}
\right.
\end{equation}
Thus, in the case of even $N$, the extremizers form a one-parameter family of
polynomials, wherein all the coefficients with odd indices are determined
uniquely, and those with even indicies depend linearly on a parameter. The formulas
for the coefficients in the problems $C_\text{max}$, $D_\text{max}$,
$D_\text{min}$ and $C_\text{min}$ are defined by 
\eqref{a-odd:case-c}-\eqref{a-even:case-c},  and
\eqref{a:case-d}-\eqref{a:case-f}, respectively.

Let us now transform the recurrence formulas, by which the coefficients of the extremizers are computed into explicit representations of the coefficients.

\subsubsection{Case a) ($A_\text{max}$, $B_\text{max}$)}
Recall
\[
y_N^{(1)}=\cos\frac{2\pi}{N+5}\quad\mbox{and set}\quad
\quad
R_N^{(1)}=\frac{8y_N^{(1)}(1-(y_N^{(1)})^2)}{N+5}.
\]
Then it follows from~\eqref{a:case-a} that $a_1=1$  and for  $k=1,\ldots,(N-1)/2$
\[
\begin{array}{l}
a_{2k+1}=(-1)^k\left(1+R_N^{(1)}\sum\limits_{j=1}^k (-1)^j U_{j-1}(y_N^{(1)})
U'_{\frac{N+3}{2}-j}(y_N^{(1)})
\right). \\
\end{array}
\]

Use Lemma~\ref{le:A8}.\textit{a)} to write
\begin{align*}
a_{2k+1}=(-1)^k&\left(1+\frac{R_N^{(1)}}{2(1-(y_N^{(1)})^2)}\right. \\
&\left.\times
\sum_{j=1}^k(-1)^j\left[\left(\frac{N+3}{2}-j\right)
U_{2j-1}(y_N^{(1)})+2U_{j-1}(y)U_j(y_N^{(1)})\right]\right).
\end{align*}
Using the identity $2U_{j-1}(x)U_j(x)=U_{2j-1}(x)+2xU^2_{j-1}(x)$, we obtain
\begin{align*}
a_{2k+1}=(-1)^k&\left(1+\frac{R_N^{(1)}}{2(1-(y_N^{(1)})^2)}\right. \\
&\left.\times\sum_{j=1}^k(-1)^j\left[\left(\frac{N+5}{2}-j\right)
U_{2j-1}(y_N^{(1)})+2y_N^{(1)}U_{j-1}^2(y_N^{(1)})\right]\right).
\end{align*}
Apply Lemma~\ref{le:A9}. By performing identical transformations, we finally obtain for $k=0,\ldots,(N-1
)/2$
\begin{align*}
a_{2k+1}=&U_{2k}(y_N^{(1)})+\frac{2}{N+5}\cdot\frac{y_N^{(1)}}{1-(y_N^{(1)})^2}\\
&\times\left(
y_N^{(1)}+U_{2k-1}(y_N^{(1)})\frac{3(y_N^{(1)})^2-1}{2(y_N^{(1)})^2}-
U_{2k}(y_N^{(1)})\frac{k(1-(y_N^{(1)})^2)+(y_N^{(1)})^2}{y_N^{(1)}}
\right).
\end{align*}
All the coefficients with even indices are equal to zero.
\subsubsection{Case d) ($B_\text{min}$)}
Recall
\[
y_N^{(2)}=\sin\frac{\pi}{N+5}\quad\mbox{and set}\quad
R_N^{(2)}=(-1)^\frac{N-1}{4}\frac{8y_N^{(2)}(1-(y_N^{(2)})^2)}{N+5}.
\]
Then, repeating the argument from the above case, we derive the analogous formula for $k=0,\ldots,(N-1
)/2.$
\begin{align*}
a_{2k+1}=&U_{2k}(y_N^{(2)})+\frac{2}{N+5}\cdot\frac{y_N^{(2)}}{1-(y_N^{(2)})^2}\\
&\times\left(
y+U_{2k-1}(y_N^{(2)})\frac{3(y_N^{(2)})^2-1}{2(y_N^{(2)})^2}-
U_{2k}(y_N^{(2)})\frac{k(1-(y_N^{(2)})^2)+(y_N^{(2)})^2}{y_N^{(2)}}
\right),
\end{align*}
All the coefficients with even indices are equal to zero.
\subsubsection{Case b) ($C_\text{max}$, $D_\text{max}$)}
Analogously to the previous argumentation recall
\[
y_N^{(3)}=\cos\frac{2\pi}{N+4} \quad\mbox{and set}\quad
R_N^{(3)}=\frac{8y_N^{(3)}(1-(y_N^{(3)})^2)}{N+4}.
\]
Then for the coefficients with odd indices, we obtain for $k=0,\ldots,(N-2
)/2$
\begin{align*}
a_{2k+1}=&U_{2k}(y_N^{(3)})+\frac{2}{N+4}\cdot\frac{y_N^{(3)}}{1-(y_N^{(3)})^2}\\
&\times\left(
y+U_{2k-1}(y_N^{(3)})\frac{3(y_N^{(3)})^2-1}{2(y_N^{(3)})^2}-
U_{2k}(y_N^{(3)})\frac{k(1-(y_N^{(3)})^2)+(y_N^{(3)})^2}{y_N^{(3)}}
\right),
\end{align*}

Compute the coefficients with even indices. By formulas~\eqref{a-even:case-c}
we get
\begin{align*}
a_{2k+2}&=(-1)^k\bigg(a_2+\frac{\tau}{2}R_N^{(3)}\bigg.\\
&\bigg.\qquad\times\sum_{j=1}^k(-1)^j
\left[U_{j-1}(y_N^{(3)})U'_{\frac{N+2}{2}-j}(y_N^{(3)})
+U_j(y)U'_{\frac{N}{2}-j}(y_N^{(3)})\right]\bigg)\\
&=(-1)^k\left(a_2+\frac{\tau}{2}R_N^{(3)}\left((-1)^j U_k(y_N^{(3)})
U'_{\frac{N}{2}-k}(y_N^{(3)})-U'_{\frac{N}{2}}(y_N^{(3)})\right)\right).    
\end{align*}
Note that 
\[
U'_\frac{N}{2}(y_N^{(3)})=\frac{y_N^{(3)}}{1-(y_N^{(3)})^2}\frac{N+4}{2},
\quad
a_2=2\tau(y_N^{(3)})^2,
\]
then for $k=0,\ldots,(N-1)/2$
\[
a_{2k+2}=\tau\frac{4y_N^{(3)}(1-(y_N^{(3)})^2)}{N+4}U_k(y_N^{(3)})
U'_{\frac{N}{2}-k}(y_N^{(3)}),
\quad\tau\in[-1,1],
\]
\subsubsection{Case f) ($D_\text{min}$)}
Recall
\[
y_N^{(4)}=\sin\frac{\pi}{N+4}\quad\mbox{and set}\quad
R_N^{(4)}=(-1)^\frac{N-2}{4}\frac{8y_N^{(4)}(1-(y_N^{(4)})^2)}{N+4}.
\]
Then it follows from~\eqref{a:case-d} that for $k=0,\ldots,(N-1)/2$
\begin{equation}
\begin{aligned}
a_{2k+1}=&U_{2k}(y_N^{(4)})+\frac{2}{N+4}\cdot\frac{y_N^{(4)}}{1-(y_N^{(4)})^2}\\
&\times\left(
y+U_{2k-1}(y_N^{(4)})\frac{3(y_N^{(4)})^2-1}{2(y_N^{(4)})^2}
-U_{2k}(y_N^{(4)})\frac{k(1-(y_N^{(4)})^2)+(y_N^{(4)})^2)}{y_N^{(4)}}
\right),
\end{aligned}
\end{equation}
\begin{equation}
a_{2k+2}=\tau\frac{4y_N^{(4)}(1-(y_N^{(4)})^2)}{N+4}
U_k(y_N^{(4)})U'_{\frac{N}{2}-k}(y_N^{(4)}),
\quad\tau\in[-1,1],
\end{equation}
\subsubsection{Cases c), e) ($A_\text{min}$, $C_\text{min}$)}
In these cases, one needs to use formulas~\eqref{a:case-e} and \eqref{a:case-f},
respectively.

The proof of Theorem~\ref{th:A} ends here.

\section{Proof of Theorem~\ref{th:B}\\
(compact representations for the extremizers)}\label{sec:4} 
The proof of Theorem~\ref{th:B} for all cases \textit{a)--e)} is conducted 
similarly and consists of three parts.

1) It is necessary to make sure that all rational functions $P(z)$ defined in 
the theorem are polynomials (of degree $N$) if they are additionally defined at
singular points by continuity (i.e., we need to show that singular points are
removable). Let $y=y_N^{(j)}$, $j=1,\ldots,6$. Note that 
\[
(1+z^4+2z^2(1-2y^2))^2=(1+z^2-2yz)^2(1+z^2+2yz)^2.
\]
The singular points of the rational function $P(z)$ will be the zeros of the
polynomials $(1+z^2\pm2yz)^2$, $(1-z^2)^k$ ($k$ equals $1$ or $3$), that is,
the set $\{\pm1,\pm e^{\pm i\alpha_N^{(j)}}\}$, where
$y=y_N^{(j)}=\cos\alpha_N^{(j)}$, with the points $\pm1$ having multiplicity 
\textit{one} in the cases \textit{a)--f)} and \textit{three} in the cases 
\textit{c), e)}, and the points $\pm e^{\pm i\alpha_N^{(j)}}$ having multiplicity
\textit{two} in all cases under consideration. Besides, there is another singular
point $\{0\}$ in the cases \textit{b), f)}. The fact that all these singular points
are removable singularities is proved in Lemmas~\ref{le:A10} and \ref{le:A12}.

2) We must ensure that these polynomials are typically real. This fact follows
from Corollary~  \ref{cor:2}.

3) Carrying out straightforward calculations, we find the third coefficient for
each constructed polynomial. In all the cases, we get $a_3=4y^2-1$ (see Theorem \ref{th:A}). Thus, Theorem~\ref{th:B} is proved.

\section{Extremal non-negative sine polynomials}\label{sec:5}

Let us show that our approach leads to the solution of the problem: What are the sine extremizers for $a_3$ which attain the bounds \eqref{eq:bounds} given by Rogosinski - Szeg\"o \cite{RS50}. Observe that on the one hand 
\begin{equation}\label{eq:pN}
   \text{Im}\{P_N(e^{it})\}=\text{Im}\{e^{it}+\sum_{j=2}^N a_je^{ijt} \}=\sin{t}+\sum_{j=2}^N a_j \sin{jt}. 
\end{equation}
On the other hand, the compact representation for $P_N(z)$ given in Theorem \ref{th:B} lead to the compact representations for $\text{Im}\{P_N(e^{it})\}.$ Introduce
\[
\hspace{-4cm}\Theta(t,y)=
\frac{2y^2(1-y^2)}{\sin t}
\left\{
\begin{array}{l}
\frac{\sin^2\frac{N+5}2 t}
{(N+5)(\cos^2t-y^2)^2},\ N\textit{ is odd},\\
\\
\frac{\sin^2\frac{N+4}2 t}
{(N+4)(\cos^2t-y^2)^2},\ N\textit{ is even},
\end{array}
\right.
\]
\bigskip
\[
\hspace{-2.2cm}\hat{\Theta}(t,y)=
\frac{2y^2(1-y^2)}{\sin^3 t}
\left\{
\begin{array}{l}
\frac{\left(\frac{N+7}2\sin\frac{N+3}2 t-\frac{N+3}2\sin\frac{N+7}2 t\right)^2}
{(N+3)(N+5)(N+7)(\cos^2 t-y^2)^2},\ N\textit{ is odd},\\
\\
\frac{\left(\frac{N+6}2\sin\frac{N+2}2 t-\frac{N+2}2\sin\frac{N+6}2 t\right)^2}
{(N+2)(N+4)(N+6)(\cos^2 t-y^2)^2},\ N\textit{ is even}.
\end{array}
\right.
\]
\\

\begin{corollary} \label{th:C}
The extremizers $P_N(z)$ characterized in Theorem \ref{th:B}, imply the following sine extremizers Im$\{P_N(e^{it})\}$  with  $ y_N^{(j)} , j=1,...,6,$ as defined in Theorem \ref{th:A}:
\\
\hspace*{.4cm} \textit{a)} Cases $A_\text{max}$, $B_\text{max}$: If \textit{$N$ is odd},  
then 
$$
\text{Im}\{P_N(e^{it})\}=\Theta(t,y_N^{(1)}).
$$

\textit{b)} Cases $C_\text{max}$, $D_\text{max}$: Let  If \textit{$N$ is odd}, then
$$\text{Im}\{P_N(e^{it})\}=\Theta(t,y_N^{(3)})(1+\tau\cos{t}),\qquad\tau\in[-1,1].
$$

\textit{c)} Case $A_\text{min}$: If \textit{$N$ and $(N+3)/2$ are odd}, then 
$$
\text{Im}\{P_N(e^{it})\}=\hat{\Theta}(t,y_N^{(5)}).
$$

\textit{d)} Case $B_\text{min}$: If \textit{$N$ is odd while $(N+3)/2$ is even}, then 
$$\text{Im}\{P_N(e^{it})\}=\Theta(t,y_N^{(2)}).
$$

\textit{e)} Case $C_\text{min}$: If \textit{$N$ is even, while $(N+2)/2$ is odd}, then
$$
\text{Im}\{P_N(e^{it})\}=\hat{\Theta}(t,y_N^{(6)})(1+\tau\cos{t}),\qquad \tau\in[-1,1].$$

\textit{f)} Case $D_\text{min}$: If \textit{$N$ and $(N+2)/2$ are even}, then 
$$
\text{Im}\{P_N(e^{it})\}=\Theta(t,y_N^{(4)})(1+\tau\cos{t}),\qquad \tau\in[-1,1].
$$
\end{corollary}

\begin{corollary} \label{cor:2}
The extremal polynomials constructed in Theorems \ref{th:A} and \ref{th:B}
 are typically real.
 \end{corollary} 

\begin{proof}
    The representations in Corollary \ref{th:C} imply that $\text{Im}\{P_N(e^{it})\}$ is non-negative $[0,\pi].$ In particular, if one takes over the notation $R_N(\cos t):= \text{Im}\{P_N(e^{it})\}/\sin t$  from Royster-Suffridge \cite[p.308]{RS70} it is clear that $R_N(\cos t)\ge 0$ on $[-\pi,\pi).$  Now, by the Main Theorem in \cite{RS70} this is equivalent to $P_N\in\mathcal{T}_N$ being typically real, thus the extremizers belong to $\mathcal T_N.$ 
\end{proof}

{
As is shown in \cite{DST25}, one can extract from the sine extremizers for the second coefficient $a_2$ new nonnegative even approximate identities. The fact that in the present situation there are essentially six sine extremizers (whereas in the $a_2$-case only two) requires a modified construction to obtain nonnegative approximate identities. This will be performed in a subsequent paper.
}

\section{Appendix}\label{sec:6}
\setcounter{theorem}{0}
\renewcommand{\thetheorem}{A.\arabic{theorem}}
\setcounter{equation}{0}
\renewcommand{\theequation}{A.\arabic{equation}}

\begin{lemma}\label{le:A3}
Denote the minimum positive root of the equation $U_j(x)=0$ by $\xi^{(j)}$ and
the minimum positive root of the equation $U'_j(x)=0$ by $\eta^{(j)}$. Then
\begin{align*}
&\eta^{(k)}<\xi^{(k-1)}<\xi^{(k)}<\eta^{(k-1)}\ \text{ when }k\textit{ is odd},\\
&\xi^{(k)}<\eta^{(k-1)}<\eta^{(k)}<\xi^{(k-1)}\ \text{ when }k\textit{ is even}.\\
\end{align*}
\end{lemma}
\begin{proof}
Let us prove the case $\eta^{(k)}<\xi^{(k-1)}$ when \textit{$k$ is odd}. Use the 
formula~\cite{DSS22,DSS23}
\[
U'_k(x)=\frac{1}{2(1-x^2)}\left((k+2)U_{k-1}(x)-kU_{k+1}(x)\right).
\]
Then 
\begin{align*}
&U'_k(0)=\frac12\left((k+2)U_{k-1}(0)-kU_{k+1}(0)\right),\\
&U'_k(\xi^{(k-1)})=\frac{-k}{2(1-(\xi^{(k-1)})^2)}U_{k+1}(\xi^{(k-1)}). 
\end{align*}
Hence, the signs of $U'_k(0)$ and $U'_k(\xi^{(k-1)})$ are opposite. This means that
the function $U'_k(x)$ has a zero in the interval $(0,\xi^{(k-1)})$, that is,
$\eta^{(k)}<\xi^{(k-1)}$. The remaining cases are treated similarly. The lemma is 
proved.
\end{proof}
As a consequence of Lemma~\ref{le:A3}, we obtain the following result.
\begin{lemma}\label{le:A4}
The maximum root $x_\text{max}$ of the equation $\det\Phi_N(x)=0$ is
\[
x_\text{max}=\left\{
\begin{array}{l}
\displaystyle
\cos\frac{2\pi}{N+5},\ N\textit{ is odd},\\
\displaystyle
\cos\frac{2\pi}{N+4},\ N\textit{ is even}.\\
\end{array}
\right.
\]
The minimum root is 
\[
x_\text{min}=\left\{
\begin{array}{ll}
\nu_1, &N\textit{ is odd},\ (N-1)/2\textit{ is odd},\\
\displaystyle
\sin\frac{\pi}{N+5},&N\textit{ is odd},\ (N-1)/2\textit{ is even},\\
\nu_2, &N\textit{ is even},\ N/2\textit{ is even},\\
\displaystyle
\sin\frac{\pi}{N+4}, &N\textit{ is even},\ N/2\textit{ is odd}.\\
\end{array}
\right.
\]
Here $\nu_1$ is the minimum positive root of the equation $U'_{(N+3)/2}(x)=0$
and $\nu_2$ is the minimum positive root of the equation $U'_{(N+2)/2}(x)=0$. 
\end{lemma}

 In the next two lemmas we denote the zero vector with  $\mathbf{0}$ and remind that the matrix $D_N(x)$ is given in \eqref{DN}.

\begin{lemma}[\cite{DSS22,DSS23}]\label{le:A5}
Let 
\[
Z^{(0)}(x)=\{\zeta_1(x),\zeta_2(x),\ldots,\zeta_N(x)\},
\quad
\zeta_k(x)=U_{k-1}(x)U_k(x),
\quad
k=1,\ldots,N;
\]
\[
x_j=\frac{\pi j}{N+2}, 
\quad
j=1,\ldots,N+1.
\]
Then
\[
D_N(x_j)\cdot Z^{(0)}(x_j)=\mathbf{0}.
\]
\end{lemma}
\begin{lemma}[\cite{DSS22,DSS23}]\label{le:A6}
Let
\[
Z^{(1)}(x)=\{\zeta_1(x),\zeta_2(x),\ldots,\zeta_N(x)\},
\]
\[
\zeta_k(x)=U_{k-1}(x)U_k(x)-U_{N-k}(x)U_{N+1-k}(x)+
\left(1-\frac{2k}{N+1}\right)U_N(x)U_{N+1}(x),
\]
\[
k=1,\ldots,N,
\]
and $\nu$ be a root of the equation $U'_{N+1}(x)=0$. Then
\[
D_N(\nu)\cdot Z^{(1)}(\nu)=\mathbf{0}.
\] 
\end{lemma}

\begin{lemma}\label{le:A7}
\textit{a)} For \textit{odd $N$}, 
\[
U_{\frac{N+1}{2}-k}\left(\cos\frac{2\pi}{N+5}\right)
=U_k\left(\cos\frac{2\pi}{N+5}\right);
\]
\textit{b)} for \textit{odd $N$} and \textit{odd $(N+1)/2$}, 
\[
U_{\frac{N+1}{2}-k}\left(\sin\frac{\pi}{N+5}\right)
=(-1)^\frac{N-1}{4} U_k\left(\sin\frac{\pi}{N+5}\right);
\]
\textit{c)} For \textit{even $N$}, 
\[
U_{\frac{N}{2}-k}\left(\cos\frac{2\pi}{N+4}\right)
=U_k\left(\cos\frac{2\pi}{N+4}\right);
\]
\textit{d)} For \textit{even $N$} and \textit{even $(N-2)/2$}, 
\[
U_{\frac{N}{2}-k}\left(\sin\frac{\pi}{N+4}\right)
=(-1)^\frac{N-2}{4}U_k\left(\sin\frac{\pi}{N+4}\right).
\]
\end{lemma}
The proof is trivial.
\begin{lemma}\label{le:A8}
\textit{a)} Let \textit{$N$ be odd}, $y=\cos2\pi/(N+5)$, then 
\[
U_{k-1}(y)U'_{\frac{N+3}{2}-k}(y)=\frac{1}{2(1-y^2)}
\left[\left(\frac{N+3}{2}-k\right)U_{2k-1}(y)+2U_{k-1}(y)U_k(y)\right];
\]
\textit{b)} let \textit{$N$ be odd, $(N+1)/2$ be odd}, $y=\sin\pi/(N+5)$, then 
\[
U_{k-1}(y)U'_{\frac{N+3}{2}-k}(y)=\frac{(-1)^\frac{N-1}{4}}{2(1-y^2)}
\left[\left(\frac{N+3}{2}-k\right)U_{2k-1}(y)+2U_{k-1}(y)U_k(y)\right];
\]
\textit{c)} let \textit{$N$ be even}, $y=\cos2\pi/(N+4)$, then 
\[
U_{k-1}(y)U'_{\frac{N+2}{2}-k}(y)=\frac{1}{2(1-y^2)}
\left[\left(\frac{N+2}{2}-k\right)U_{2k-1}(y)+2U_{k-1}(y)U_k(y)\right];
\]
\textit{d)} let \textit{$N$ be even, $(N-2)/2$ be even}, $y=\sin\pi/(N+4)$, then 
\[
U_{k-1}(y)U'_{\frac{N+2}{2}-k}(y)=\frac{(-1)^\frac{N-2}{4}}{2(1-y^2)}
\left[\left(\frac{N+2}{2}-k\right)U_{2k-1}(y)+2U_{k-1}(y)U_k(y)\right].
\]
\end{lemma}
\begin{proof}
\textit{a)} Apply the formula
\[
U'_k(x)=\frac{1}{2(1-x^2)}((k+2)U_{k-1}(x)-kU_{k+1}(x)).
\]
Then 
\begin{align*}
U_{k-1}(y)&U'_{\frac{N+3}{2}-k}(y) \\
&=\frac{U_{k-1}(y)}{2(1-y^2)}
\left[\left(\frac{N+3}{2}-k\right)\left(U_{\frac{N+1}{2}-k}(y)
-U_{\frac{N+3}{2}-k}(y)\right)+2U_{\frac{N+1}{2}-k}(y)\right].
\end{align*}
By Lemma~\ref{le:A7} it follows that
\[
U_{k-1}(y)U'_{\frac{N+3}{2}-k}(y)=\frac{U_{k-1}(y)}{2(1-y^2)}
\left[\left(\frac{N+3}{2}-k\right)(U_k(y)-U_{k-2}(y))+2U_k(y)\right].
\]
By using the obvious identity $U_{k-1}(x)U_k(x)-U_{k-1}(x)U_{k-2}(x)=U_{2k-1}(x)$,
we obtain the statement \textit{a)} of Lemma~\ref{le:A8}. The rest of the 
statements are proved in the same way. 
\end{proof}
\begin{lemma}\label{le:A9}
The following identities hold:

\textit{a)} 
$
\displaystyle
\sum_{j=1}^k (-1)^j U_{2j-1}(x)=\frac{(-1)^k}{2x}(U_{2k}(x)-(-1)^k);
$

\textit{b)} 
$
\displaystyle
\sum_{j=1}^k (-1)^j jU_{2j-1}(x)=\frac{(-1)^k}{4(1-x^2)}
((k+1)U_{2k-1}(x)+kU_{2k+1}(x));
$

\textit{c)} 
$
\displaystyle
\sum_{j=1}^k (-1)^j U_{j-1}^2(x)=\frac{(-1)^k}{2x}U_{k-1}(x)U_k(x).
$
\end{lemma}
\begin{proof}[Proof\nopunct]
follows from~\cite[Section~4.4.1, formulas 6, 8; Section 4.4.2, formula 3]{PBM86}.
\end{proof}
\begin{lemma}\label{le:A10}
Let
\[
G_1(z,y)=\left\{
\begin{array}{l}
\displaystyle
\frac{32}{N+5}y^2(1-y^2)\frac{z^5}{1-z^2}
\frac{1-z^{N+5}}{(1+z^4+2(1-2y^2)z^2)^2},\ N\textit{ is odd},\\
\displaystyle
\frac{32}{N+4}y^2(1-y^2)\frac{z^5}{1-z^2}
\frac{1-z^{N+4}}{(1+z^4+2(1-2y^2)z^2)^2},\ N\textit{ is even},
\end{array}
\right.
\]
\[
G_2(z,y)=\frac{z+z^3}{1+z^4+2(1-2y^2)z^2}.
\]
Consider the following four cases:
\begin{enumerate}
\item[\textit{i)}] \textit{$N$ is odd}, $y_N=\cos2\pi/(N+5)$;
\item[\textit{ii)}] \textit{$N$ is odd, $s=(N+3)/2$ is even}, 
$y_N=\cos\pi s/(N+5)=\sin\pi/(N+5)$;
\item[\textit{iii)}] \textit{$N$ is even}, $y_N=\cos2\pi/(N+4)$;
\item[\textit{iv)}] \textit{$N$ is even, $s=(N+2)/2$ is even}, 
$y_N=\cos\pi s/(N+4)=\sin\pi/(N+4)$.
\end{enumerate}
Then in cases \textit{i), ii)}, the rational function 
$P(z)=G_1(z,y_N)+G_2(z,y_N)$ has only removable singularities at the points
$\{\pm1,\pm e^{\pm i\arccos(y_N)}\}$; in the cases \textit{iii), iv)}, the function
\begin{equation}\label{eq:A10}
P(z)=(G_1(z,y_N)+G_2(z,y_N))\left(1+\tau\frac{z+z^{-1}}{2}\right)-\frac{\tau}{2}
\end{equation}
has as singular points only removable singularities 
$\{0,\pm1,\pm e^{\pm i\arccos(y_N)}\}$.
\end{lemma}
\begin{proof}
The fact that the singular points $0,\pm1$ are removable singularities is verified
easily. Thus there remains to show that the point $z^*=e^{i\arccos(y_N)}$ is a removable
singularity of the function~\eqref{eq:A10} in the case \textit{iv)}. Note that all the
other cases are treated similarly. We need to prove that there exists a
finite limit $\lim_{z\to z^*}P(z)$, while without loss of generality we may assume
that $\tau=0$. To do this, let us introduce the auxiliary function 
\[
V(z)=z^{-4}(1+z^2)(1-z^2)(1+z^2+2yz)(1+z^2-2yz)
+\frac{32}{N+4}y^2(1-y^2)(1-z^{N+4}),
\]
$y=y_N$, and show that $V(z^*)=0$, $V'(z^*)=0$. 

Since $s=(N+2)/2$ is even, then $(z^*)^{N+4}=e^{i\pi s}=1$. Calculate
\begin{align*}
V'(z)&=[z^{-4}(1+z^2)(1-z^2)(1+z^2+2yz)]'(1+z^2-2yz)\\
&\hspace{1cm} \quad+2(z-y)z^{-4}(1+z^2)(1-z^2)(1+z^2+2yz)\\
&\hspace{2.2cm}  \quad-\frac{32}{N+4}y^2(1-y^2)(N+4)z^{-1}z^{N+4}.
\end{align*}
Then
\[
V'(z^*)=\frac{1}{z^*}
\left[
-32y^2(1-y^2)+2(y+i\sqrt{1-y^2}-y)2y(-2i)\sqrt{1-y^2}(2y+2y)
\right]=0.
\]
The lemma is proved.
\end{proof}
We will need the following auxiliary lemma.
\begin{lemma}\label{le:A11}
Let
\[
b\sin{at}-a\sin{bt}=0,
\]
then
\begin{enumerate}
\item[\textit{i)}] 
$b^2(1-\cos{2at})+a^2(1-\cos{2bt})-2ab(\cos(a-b)t-\cos(a+b)t)=0$,
\item[\textit{ii)}] $b^2\sin{2at}+a^2\sin{2bt}-2ab\sin(a+b)t=0$,
\item[\textit{iii)}] 
$b\sin{2at}+a\sin{2bt}-(a+b)\sin(a+b)t+(b-a)\sin(b-a)t=0$,
\item[\textit{iv)}] 
$b\cos{2at}+a\cos{2bt}-(a+b)\cos(a+b)t+(a+b)(-1+\cos(a-b)t)=0$.
\end{enumerate}
\end{lemma}
\begin{proof}[Proof\nopunct]
follows from the identities
\[
b^2(1-\cos{2at})+a^2(1-\cos{2bt})-2ab(\cos(b-a)t-\cos(a+b)t)
=2(b\sin{at}-a\sin{bt})^2,
\]
\[
b^2\sin{2at}+a^2\sin{2bt}-2ab\sin(a+b)t
=2(b\sin{at}-a\sin{bt})(b\cos{at}-a\cos{bt}),
\]
\begin{gather*}
b\sin{2at}+a\sin{2bt}-(a+b)\sin(a+b)t+(b-a)\sin(b-a)t\\
=2(b\sin{at}-a\sin{bt})(\cos{at}-\cos{bt})
\end{gather*}
and 
\begin{gather*}
b\cos{2at}+a\cos{2bt}-(a+b)\cos(a+b)t+(a+b)(-1+\cos(a-b)t)\\
=-2(b\sin{at}-a\sin{bt})(\sin{at}-\sin{bt}).
\end{gather*}
\end{proof}
\begin{lemma}\label{le:A12}
Let 
\[
G_3(z,y)=\left\{
\begin{array}{l}
-\frac{32y^2(1-y^2)}{(N+3)(N+5)(N+7)}\frac{z^7}{(1-z^2)^3}\\
\quad\times\frac{(N+7)^2(1-z^{N+3})+(N+3)^2(1-z^{N+7})-2(N+3)(N+7)(z^2-z^{N+5})}
{(1+z^4+2(1-2y^2)z^2)^2},\ N\textit{odd},\\
-\frac{32y^2(1-y^2)}{(N+2)(N+4)(N+6)}\frac{z^7}{(1-z^2)^3}\\
\quad\times\frac{(N+6)^2(1-z^{N+2})+(N+2)^2(1-z^{N+6})-2(N+2)(N+6)(z^2-z^{N+4})}
{(1+z^4+2(1-2y^2)z^2)^2},\ N\textit{even},
\end{array}
\right.
\]
\[
G_4(z,y)=\frac{(z+z^3)(1+z^8-\gamma_1(y)(z^2+z^6)+\gamma_2(y)z^4)}
{(1-z^2)^2(1+z^4+2(1-2y^2)z^2)^2},
\]
where 
\[
\gamma_1(y)=4y^2,
\quad
\gamma_2(y)=\left\{
\begin{array}{l}
2\left(-\frac{16}{N+5}y^4+4\left(1+\frac{4}{N+5}\right)y^2-1\right),
\ N\textit{odd},\\
2\left(-\frac{16}{N+4}y^4+4\left(1+\frac{4}{N+4}\right)y^2-1\right),
\ N\textit{even}.\\
\end{array}
\right.
\]
Consider the two cases:
\begin{enumerate}
\item[\textit{i)}] \textit{$N$ is odd, $(N+3)/2$ is odd}, $y_N=\cos\alpha_N$,
$U'_{(N+3)/2}(y_N)=0$;
\item[\textit{ii)}] \textit{$N$ is even, $(N+2)/2$ is odd}, $y_N=\cos\alpha_N$,
$U'_{(N+2)/2}(y_N)=0$. 
\end{enumerate}

Then in the case \textit{i)}, the polynomial $P(z)=G_3(z,y_N)+G_4(z,y_N)$ has as 
singular points only removable singularities $\{\pm1,\pm e^{\pm i\alpha_N}\}$;
in the case \textit{ii)}, the polynomial
\[
P(z)=(G_3(z,y_N)+G_4(z,y_N))\left(1+\tau\frac{z+z^{-1}}{2}\right)-\frac{\tau}{2}
\]
has only removable singularities at $\{0,\pm1,\pm e^{\pm i\alpha_N}\}$.
\end{lemma}
\begin{proof}
The singilar points $0,\pm1$ are removable singularities of the rational function
$G_3(z,y_N)+G_4(z,y_N)$, which is checked by direct calculations. Let us show that
the point $z^*=e^{i\alpha}$ is a removable singularity of the polynomial 
$P(z)=G_3(z,y)+G_4(z,y)$ ($\alpha=\alpha_N$, $y=y_N$, case~\textit{i)}). The
parameter $\alpha$ satisfies the equation 
\[
(N+7)\sin\frac{N+3}{2}\alpha-(N+3)\sin\frac{N+7}{2}\alpha=0.
\]
Note that case~\textit{ii)} reduces to the previous case. We need to show that 
there is a finite limit $\lim_{z\to z^*}P(z)$. To do this, we introduce the 
auxiliary function 
\begin{align*}
W(z)=z^{-6}(1-z^4)(1+&z^8-\gamma_1(z^2+z^6)+\gamma_2z^4)
-Q\left[(N+7)^2(1-z^{N+3})\right.\\
&\left.+(N+3)^2(1-z^{N+7})-2(N+3)(N+7)(z^2-z^{N+5})\right],
\end{align*}
where 
\[
\gamma_1=4y^2,
\quad
\gamma_2=2\left(-\frac{16}{N+5}y^4+4\left(1+\frac{4}{N+5}\right)y^2-1\right),
\]
\[
Q=\frac{32y^2(1-y^2)}{(N+3)(N+5)(N+7)},
\]
and show that $W(z^*)=0$, $W'(z^*)=0$.

Set
\begin{align*}
&W_1(z)=z^{-6}(1-z^4)(1+z^8-\gamma_1(z^2+z^6)+\gamma_2z^4),\\
&\begin{aligned}
W_2(z)=Q\left[(N+7)^2(1-z^{N+3})\right.+&(N+3)^2(1-z^{N+7})\\
&\left.-2(N+3)(N+7)(z^2-z^{N+5})\right].
\end{aligned}
\end{align*}
Compute $W_1(e^{i\alpha})=-2i\sin2\alpha(2\cos4\alpha-
2\gamma_1\cos2\alpha+\gamma_2)$, whence we obtain, in view of $y=\cos\alpha$, that
\[
W_1(e^{i\alpha})=-i\frac{2^7}{N+5}y^3(1-y^2)\sqrt{1-y^2}.
\]
Further,
\begin{align*}
W_2(e^{i\alpha})&=Q(N+7)^2(1-\cos(N+3)\alpha)+(N+3)^2(1-\cos(N+7)\alpha)\\
&\quad+2(N+3)(N+7)(\cos2\alpha-\cos(N+5)\alpha)-
iQ\left[(N+7)^2\sin(N+3)\alpha\right.\\
&\quad\left.+(N+3)^2\sin(N+7)\alpha+2(N+3)(N+7)(\sin2\alpha
-\sin(N+5)\alpha)\right]\\
&=\frac{2^7}{N+5}iy^3(1-y^2)\sqrt{1-y^2}.
\end{align*}
In the last transformation, we applied Lemma~\ref{le:A11}, 
formulas~\textit{i), ii)} with $a=(N+3)/2$, $b=(N+7)/2$, $y=\cos\alpha$. 
Thus, it follows that $W(e^{i\alpha})=0$. 

Now we show that $W'(e^{i\alpha})=0$. We consecutively calculate 
\begin{align*}
&\begin{aligned}
W'_1(z)&=-\frac{2}{z}(z^{-2}+z^2)(z^{-4}+z^4-\gamma_1(z^{-2}+z^2)+\gamma_2)\\
&\quad+\frac{1}{z}(z^{-2}-z^2)(-4(z^{-4}-z^4)+2\gamma_1(z^{-2}-z^2)),
\end{aligned}\\
&\begin{aligned}
W'_1(e^{i\alpha})&=-4e^{-i\alpha}\cos2\alpha(2\cos4\alpha
-2\gamma_1\cos2\alpha+\gamma_2)\\
&\quad-2ie^{-i\alpha}\sin{2\alpha}(-8i\sin{4\alpha}-4i\gamma_1\sin{2\alpha}).
\end{aligned}
\end{align*}
Taking into account the relation $y=\cos\alpha$, we get
\[
W'_1(e^{i\alpha})=-\frac{2^7}{N+5}e^{-i\alpha}y^2(1-y^2)(1-2y^2)
-2^7e^{-i\alpha}y^2(1-y^2)^2,
\]
and further,
\[
W'_2(z)=\frac{1}{z}\frac{2^5}{N+5}y^2(1-y^2)
\left[(N+7)z^{N+3}+(N+3)z^{N+7}+4z^2-2(N+5)z^{N+5}\right].
\]
An application of Lemma~\ref{le:A11}, formulas~\textit{iii), iv)} with $a=(N+3)/2$ and 
$b=(N+7)/2$, yields 
\[
W'_2(e^{i\alpha})=e^{-i\alpha}\left(-\frac{2^7}{N+5}y^2(1-y^2)(1-2y^2)
+2^7y^2(1-y^2)^2\right).
\]
Finally, this results in $W'_1(e^{i\alpha})+W'_2(e^{i\alpha})=0$. The lemma is
proved.

\end{proof}
\bibliographystyle{abbrv}
\bibliography{bibliography}
\end{document}